\documentclass[moor,nonblindrev]{informs1mod} 

\usepackage{tablefootnote}
\usepackage{float}
\usepackage[dvipsnames]{xcolor}
\usepackage{mathrsfs}	
\usepackage[shortlabels,inline]{enumitem}
\setlist{noitemsep}
\usepackage{bbm}
\usepackage{stmaryrd}
\usepackage{dsfont}
\usepackage{anyfontsize}
\newlist{glossa}{description}{1}
\newlist{glossb}{description}{1}
\setlist[glossa, 1]{style = standard, labelwidth=2.5cm, leftmargin = 2.5cm, labelsep = 0pt, itemsep = 1pt}
\setlist[glossb, 1]{style = standard, labelwidth=4.5cm,   leftmargin = 4.5cm  , labelsep = 0pt, itemsep = 1pt}
\counterwithin*{section}{part}

\usepackage{natbib}
 \NatBibNumeric
 \def\bibfont{\small}%
 \bibpunct[, ]{[}{]}{,}{n}{}{,}%

\usepackage[colorlinks=true,breaklinks=true,urlcolor=blue,
     citecolor=blue,linkcolor=blue,draft=false]{hyperref}
\usepackage{bookmark}		
\bookmarksetup{open,numbered}

\def\EMAIL#1{\href{mailto:#1}{#1}}

\TheoremsNumberedBySection  

\EquationsNumberedBySection 

\makeatletter
\newcommand*{\bigcdot}{}%
\DeclareRobustCommand*{\bigcdot}{%
  \mathbin{\mathpalette\bigcdot@{}}%
}
\newcommand*{\bigcdot@scalefactor}{.7}
\newcommand*{\bigcdot@widthfactor}{1.15}
\newcommand*{\bigcdot@}[2]{%
  \sbox0{$#1\vcenter{}$}%
  \sbox2{$#1\cdot\m@th$}%
  \hbox to \bigcdot@widthfactor\wd2{%
    \hfil
    \raise\ht0\hbox{%
      \scalebox{\bigcdot@scalefactor}{%
        \lower\ht0\hbox{$#1\bullet\m@th$}%
      }%
    }%
    \hfil
  }%
}
\makeatother

\makeatletter
\newcommand*\rel@kern[1]{\kern#1\dimexpr\macc@kerna}
\newcommand*\widebar[1]{%
  \begingroup
  \def\mathaccent##1##2{%
    \rel@kern{0.8}%
    \overline{\rel@kern{-0.8}\macc@nucleus\rel@kern{0.2}}%
    \rel@kern{-0.2}%
  }%
  \macc@depth\@ne
  \let\math@bgroup\@empty \let\math@egroup\macc@set@skewchar
  \mathsurround\z@ \frozen@everymath{\mathgroup\macc@group\relax}%
  \macc@set@skewchar\relax
  \let\mathaccentV\macc@nested@a
  \macc@nested@a\relax111{#1}%
  \endgroup
}
\makeatother

\DeclareMathAlphabet {\scr }{U}{scr}{m}{n}
\DeclareFontFamily {U}{scr}{}
\DeclareFontShape {U}{scr}{m}{n}{<-6>rsfs5<6-8>rsfs7<8->rsfs10}{}
\DeclareFontShape {U}{scr}{b}{n}{<-6>rsfs5<6-8>rsfs7<8->rsfs10}{}

\newcommand{\bpf}{\noindent {\em Proof.}\ }
\newcommand{\ep}{\hfill $\square $\vspace{0.2ex}}
\renewcommand{\qed}{ $\hfill\square $}
\newcommand{\Exp}{\mathscr{E}}

\newcommand{\sint}{\bigcdot}
\newcommand{\auf}{[\![}
\newcommand{\zu}{]\!]}
\newcommand{\Log}{\mathcal{L}}

\newcommand{\sigalg}[1]{\mathscr{#1}}
\newcommand{\CK}{\mathrm{CK}}
\renewcommand{\P}{\mathsf{P}}
\newcommand{\E}{\textsf{\upshape E}}
\renewcommand{\d}{\mathrm{d}}

\DeclareMathOperator{\esssup}{ess\ sup}
\DeclareMathOperator{\essinf}{ess\ inf}
\DeclareMathOperator{\Var}{\textsf{\upshape Var}}
\DeclareMathOperator{\diag}{diag}

\newcommand{\tame}{\Theta}
\newcommand{\adm}{\smash{\mkern1mu\widebar{\mkern-1mu\tame\mkern-1mu}\mkern1mu}}
\newcommand{\vt}{\vartheta}
\newcommand{\vp}{\varphi}

\newcommand{\eps}{\varepsilon}

\newcommand{\veps}{\epsilon}

\newcommand{\F}{\sigalg{F}}
\newcommand{\R}{\mathbb{R}}
\newcommand{\N}{\mathbb{N}}
\newcommand{\Cx}{\mathbb{C}}
\newcommand{\cM}{\mathcal{M}}
\newcommand{\T}{\top}
\newcommand{\hY}{\hat{Y}}
\newcommand{\hZ}{\hat{Z}}
\newcommand{\hS}{\hat{S}}
\newcommand{\hP}{\hat{\P}}
\newcommand{\hPstar}{\hat{\P}^\star}
\newcommand{\tP}{\tilde{\P}}
\newcommand{\hV}{\hat{V}}
\newcommand{\hv}{\hat{v}}
\newcommand{\hc}{\hat{c}}
\newcommand{\hb}{\hat{b}}
\newcommand{\ha}{\hat{a}}

\newcommand{\hL}{\hat{L}}
\newcommand{\hH}{\hat{H}}
\newcommand{\hxi}{\hat{\xi}}
\newcommand{\hvp}{\hat{\vp}}
\newcommand{\hvt}{\hat{\vt}}
\newcommand{\hXi}{\widehat{\Xi}}
\newcommand{\hG}{\hat{G}}
\newcommand{\hW}{\hat{W}}
\newcommand{\e}{\mathrm{e}}
\renewcommand{\dag}{{-1}}
\newcommand{\rank}{r}
\newcommand{\X}{\mathcal{X}}
\newcommand{\Y}{\mathcal{Y}}
\newcommand{\bX}{\widebar{\mathcal{X}}}
\newcommand{\bY}{\widebar{\mathcal{Y}}}
\newcommand{\A}{\mathcal{A}}
\DeclareMathOperator{\Ran}{\mathcal{R}}
\DeclareMathOperator{\Null}{\mathcal{N}}
\newcommand{\NullStrategy}{\mathscr{N}\!}
\newcommand{\ones}{\mathbbm{1}}
\newcommand{\one}{\mathbf{1}}
\def\MR#1{\href{http://www.ams.org/mathscinet-getitem?mr=#1}{MR-#1}}
\def\ARXIV#1{\href{https://arxiv.org/abs/#1}{arXiv:#1}}
\def\DOI#1{\href{https://doi.org/#1}{doi:#1}}

\begin{document}

\RUNAUTHOR{\v{C}ern\'{y}, Czichowsky, and Kallsen}

\RUNTITLE{Numeraire-invariant quadratic hedging and mean--variance portfolio allocation}

\TITLE{Numeraire-invariant quadratic hedging and \\ mean--variance portfolio allocation}

\ARTICLEAUTHORS{%
\AUTHOR{Ale\v{s} \v{C}ern\'{y}}
\AFF{Bayes Business School, City, University of London, London EC1Y 8TZ, United Kingdom, \EMAIL{ales.cerny@city.ac.uk}}
\AUTHOR{Christoph Czichowsky}
\AFF{London School of Economics and Political Science, London WC2A 2AE, United Kingdom,  \EMAIL{c.czichowsky@lse.ac.uk}}
\AUTHOR{Jan Kallsen}
\AFF{Christian-Albrechts-Universit\"{a}t zu Kiel, 24118 Kiel, Germany, \EMAIL{kallsen@math.uni-kiel.de}}
} 

\ABSTRACT{%
The paper investigates quadratic hedging in a semimartingale market that does not necessarily contain a risk-free asset. An equivalence result for hedging with and without numeraire change is established. This permits direct computation of the optimal strategy without choosing a reference asset and/or performing a numeraire change. New explicit expressions for optimal strategies are obtained, featuring the use of oblique projections that provide unified treatment of the case with and without a risk-free asset. The analysis yields a streamlined computation of the efficient frontier for the pure investment problem in terms of three easily interpreted processes. The main result advances our understanding of the efficient frontier formation in the most general case where a risk-free asset may not be present. Several illustrations of the numeraire-invariant approach are given.
}%

\MSCCLASS{Primary: 91G10; 93E20; secondary 15A10\smallskip}
\ORMSCLASS{Primary: Finance/portfolio; dynamic programming/optimal control/applications; secondary: mathematics/matrices}
\HISTORY{Author accepted manuscript, April 7, 2023; final published version \DOI{10.1287/moor.2023.1374}, \emph{Mathematics of Operations Research}, 2024, 49(2), 752-781; post-publication edit in Proposition~\ref{P:new2} and Theorem~\ref{T:constrained q} (weaker assumptions, stronger statement of \eqref{eq: modifP}) and updated references, July 7, 2025.}

\maketitle

\setlength{\abovedisplayskip}{9pt}
\setlength{\belowdisplayskip}{9pt}

\section{Introduction.}
Markowitz's  mean--variance portfolio optimization is one of the pillars of modern financial economic theory, underpinning large
parts of investment practice. While this problem has been considered in one period with or without a risk-free asset, a multi-period or continuous-time analysis usually involves the existence of a risk-free asset that is taken as the numeraire and/or as the reference asset. The choice of numeraire then typically affects both the problem formulation via the set of admissible trading strategies and the formulae for optimal trading strategies through the input quantities they are derived from.

In this paper, we consider  mean--variance portfolio optimization and quadratic hedging without necessarily assuming a risk-free asset or choosing a numeraire asset. Our main contributions are a numeraire-invariant problem formulation including a symmetric definition of admissible trading strategies (Definition~\ref{D:adm}),
sufficient conditions for the existence of optimal trading strategies (Theorem~\ref{T:closedness}), expressions for the optimal trading strategies that do not require the choice of a reference asset and/or numeraire change (Theorems~\ref{T:main} and~\ref{T:explicit}), and an equivalence result for hedging with and without numeraire change (see (\ref{eq:3}--\ref{eq:4}) and Proposition \ref{P:ThetaEqual}). Our results apply in both multi-period discrete time and continuous time settings and allow computing optimal strategies in concrete models.
In particular, the framework yields the optimal strategy for  mean--variance portfolio optimization and quadratic hedging in the case of only risky assets.

Our analysis provides a streamlined computation of the efficient frontier for the pure investment problem in terms of three easily interpreted processes: the \emph{opportunity process} $L$ that measures the smallest second moment of a fully invested portfolio, a \emph{tracking process} $V(1)$ used for hedging the constant payoff 1, and the corresponding \emph{minimal expected squared hedging error} $\eps^2(1)$. The efficient frontier is characterized by the equation
\begin{equation}\label{eq:efintro}
\Var(R) =\frac{L_0\eps_0^2(1)}{L_0V^2_0(1)+\eps_0^2(1)}+\left(\frac{1}{1-L_0V^2_0(1)-\eps_0^2(1)}-1\right)\left( \E[R]-\frac{L_0V_0(1)}{L_0V^2_0(1)+\eps_0^2(1)}\right) ^{2},
\end{equation}
linking the variance of any efficient payoff $R$ to its expected value; see Section~\ref{S:efficient frontier}.
Theorems~\ref{T:main} and~\ref{T:explicit} describe how the three main ingredients are obtained in a semimartingale model. Section~\ref{S:examples} provides numerical illustrations of the streamlined approach.

The novel explicit formulae in Theorems~\ref{T:main} and~\ref{T:explicit} reveal  the key role that certain \emph{oblique projectors} play in the solution to mean--variance hedging and portfolio optimization. 
They apply universally in discrete as well as continuous time and, formulated in undiscounted terms, allow for a unified treatment of the cases with and without a risk-free asset. As an easy corollary we are able to obtain the generalization of the models in Li and Ng \cite{li.ng.00} and Yao et al. \cite{yao.al.14} to asset returns driven by arbitrary square-integrable processes with independent increments, see Section~\ref{S:examples}.

The remainder of the paper is organized as follows. In the rest of this introduction we outline the main research questions arising from  numeraire change in quadratic hedging. After establishing notation in Section~\ref{S:notation}, Section~\ref{S:adm} introduces a symmetric definition of admissibility and studies its consequences.  
In particular, Theorem~\ref{T:main X=1} in Subsection~\ref{SS:X=1} and Corollary~\ref{C:mainhP} in Subsection~\ref{SS:hX=1} give symmetric versions of the classical optimal hedging results with constant risk-free asset before and after numeraire change, respectively. Section~\ref{S:main} is devoted to optimal hedging without numeraire change; here Theorems~\ref{T:main} and~\ref{T:explicit} contain the main results of the paper. The short Section~\ref{S:efficient frontier} summarizes known facts about efficient frontiers, paving the way for Section~\ref{S:examples} that concludes with three illustrative examples of the numeraire-invariant approach without discounting. Appendix~\ref{App:A} outlines the connection among affinely constrained quadratic optimization, oblique projections, and pseudoinverses; Appendix~\ref{S:proof X=1} proves the auxiliary Theorem~\ref{T:main X=1}; Appendices~\ref{S:proof main} and~\ref{App:D} contain the proof and auxiliary statements for the main Theorem~\ref{T:main}.

\subsection{Numeraire invariance.}
We shall now outline the main building blocks of our approach. Let $(\Omega,\F,\left\{ \F_{t}\right\} _{t\in \lbrack 0,T]},\P)$ be a filtered probability space with $\F_0$ trivial. 
Our aim is to study the quadratic hedging problem
\begin{equation}\label{eq:1new}
\min_{\vt\in\adm_v(S,\P)}\E[(\vt_T S_T-H)^2], 
\end{equation}
where $H\in L^2(\P)$ is the payoff to be hedged, $S$ an $\R^d$-valued price process of $d\geq 2$ traded assets, $v\in\R$ the initial capital, and $\adm_v(S,\P)$ some set of self-financing trading strategies with initial value $v\in\R$.
Recall that with Harrison and Kreps \cite{harrison.kreps.79}, an $\R^d$-valued strategy $\vt$ that is integrable with respect to $S$ is called self-financing if its value satisfies
\begin{equation}\label{eq:SF1}
 \vt_t S_t = \vt_0S_0 + \vt \sint S_t,\quad 0\leq t\leq T,
\end{equation}
where  $\vt_t S_t =\sum_{i=1}^d \vt_t^i S_t^i$ and $\vt\sint S_t$ stands for the stochastic integral $\int_{(0,t]}\vt_u \d S_u$.

In the literature, \eqref{eq:1new} is typically considered only in discounted form. Explicitly or implicitly, it is assumed that
\begin{equation}\label{eq:S}
S = (X,Y)
\end{equation}
with $(d-1)$-dimensional $Y$ and constant $X=1$. Writing self-financing strategies in the form
$\vt=(\chi,\eta)$ with $\chi=\vt^1$ and $\eta=\vt^{2:d}=(\vt^2,\dots,\vt^d)$,
the value process \eqref{eq:SF1} of $\vt$ can in this case be rewritten as
\[\vt S=v+\eta\sint Y,\]
where $v=\vt_0 S_0$ denotes the initial wealth. Problem \eqref{eq:1new} then turns into
\begin{equation}\label{eq:1}
\min_{\eta \in \adm^{\CK}(Y,\P)}\E[(v+\eta \sint Y_{T}-H)^{2}],
\end{equation}
where $\adm^{\CK}(Y,\P)$ denotes some appropriately chosen set of admissible trading strategies (see, e.g., \cite{cerny.kallsen.07}). The reader is referred to Schweizer \cite{schweizer.10} for a recent overview of quadratic hedging.

The aim of this paper is to consider more general $S=(X,Y)$, where we only require that the price of the first asset $X$ and its left limit  $X_-$ are strictly positive.  To make progress in this setting, we use $X$ as a \emph{numeraire} and denote the discounted quantities with a hat $\hat{\textcolor{white}{v}}$, writing
\[\hv=\frac{v}{X_0},\quad \hS=\frac{S}{X},\quad\hY=\frac{Y}{X},\quad \text{etc.}\]
The discounted form of the self-financing condition now reads
\begin{equation}\label{eq:SF2}
 \vt \hS = \vt_0\hS_0 + \vt \sint \hS = \hv +\eta\sint \hY.
\end{equation}

In order to examine the link between the discounted and the undiscounted version,
assume $X_{T}\in L^2(\P)$ and let $\hP$ denote a new absolutely continuous probability measure defined by 
\begin{align}\label{eq:hP}
\frac{\d\hP}{\d\P} = \frac{X^2_{T}}{\E[X_{T}^{2}]}.
\end{align}
Observe that $\hH= \frac{H}{X_T}$ is in $L^2(\hP)$ if and only if $H\in L^2(\P)$. Furthermore, $\vt_T S_{T}$ is in $L^{2}(\P)$ if and only if $\vt_T \hS_{T}\in L^{2}(\hP)$. This allows us to write
\begin{equation}\label{eq:MeasChg2}
\E[(\vt_T S_{T}-H)^{2}]=\E[X_{T}^{2}]\E^{\hP}[(\vt_T \hS_{T}-\hH)^{2}],
\end{equation}
whenever one of the two expressions is well-defined. Let us now denote by $\adm_v(S,\P)$ the set of admissible quadratic hedging strategies for asset $S$ with initial capital $v$ under measure $\P$. The equality \eqref{eq:MeasChg2} now yields equivalence of hedging with and without change of numeraire, namely
\begin{equation}\label{eq:3}
\underbrace{\min_{\vt \in \adm_{v}(S,\P)} \E[(\vt_{T}S_{T}-H)^{2}]}_{\text{undiscounted problem $(S,H,\P)$}} 
=\E[X_{T}^{2}]
\underbrace{\min_{\vt \in \adm_{\hv }(\hS,\hP)}\E^{\hP}[(\vt_T \hS_{T}-\hH)^{2}]}
_{\text{discounted problem $(\hS,\hH,\hP)$}}
\end{equation}
and
\begin{equation}\label{eq:3.5}
\underbrace{\argmin_{\vt \in \adm_{v}(S,\P)} \E[(\vt_{T}S_{T}-H)^{2}]}_{\text{optimal strategy for $(S,H,\P)$}} 
=
\underbrace{\argmin_{\vt \in \adm_{\hv }(\hS,\hP)}\E^{\hP}[(\vt_T \hS_{T}-\hH)^{2}]}
_{\text{optimal strategy for $(\hS,\hH,\hP)$}},
\end{equation}
provided that
\begin{equation}\label{eq:4}
\adm_{v}(S,\P)=\adm_{\hv }(\hS,\hP)
\end{equation}
holds. This crucial last relation \eqref{eq:4} is shown in Proposition \ref{P:ThetaEqual} in our setup.

The solution of the second optimization in \eqref{eq:3.5} is known by virtue of the results in \v{C}ern\'y and Kallsen \cite{cerny.kallsen.07}, where thanks to \eqref{eq:SF2} one can naturally identify $\adm_{\hv }(\hS,\hP)$ with $\adm^\CK(\hY,\hP)$. Therefore, \emph{if} one decides to be guided solely by the discounted optimization in \eqref{eq:3.5}, one ought to use \eqref{eq:4} as the \emph{definition} of what is admissible in the undiscounted problem. This is certainly possible if $X$ is fixed \emph{a priori}. But if one wants to make a more fundamental progress, such approach is problematic on two grounds: (i) for some choices of the numeraire $X$ the discounted price process $\hS$ may not be locally square-integrable under $\hP$, making the discounted optimization undefined; (ii) it is in principle possible that two different choices of $X$ yield two \emph{different} sets of $\adm_{v}(S,\P)$ via \eqref{eq:4}.  

In this paper, we opt to give a direct meaning to $\adm_{v}(S,\P)$. The theory of quadratic hedging (and more broadly the utility maximization literature, e.g., \cite{biagini.cerny.20}) is concerned with designing a set of admissible strategies that precludes creation of wealth out of nothing but is rich enough to contain the optimizer. By extending the techniques of \cite{cerny.kallsen.07} we provide a \emph{symmetric} definition (Subsection~\ref{SS:sf} and Definition~\ref{D:adm}) of the set $\adm_{v}(S,\P)$ as an appropriately specified $L^{2}(\P)$-closure of highly integrable\ $d$-dimensional self-financing strategies. This guarantees, under suitable no-arbitrage conditions, that the set of admissible terminal wealths is $L^{2}(\P)$-closed and therefore an optimizer exists for the undiscounted problem (Theorem~\ref{T:closedness}).

Although the undiscounted problem now has a solution thanks to Theorem~\ref{T:closedness}, the only way to compute the solution is by passing to the 
$(\hS,\hH,\hP)$ formulation in \eqref{eq:3.5} because, at present, there is no theory that can deal with the undiscounted problem directly.
To gain access to the discounted problem \eqref{eq:3.5}, the required equality of strategies in \eqref{eq:4} calls for certain restrictions on the numeraire $X$, which we formalize in the concept of the  `nice numeraire' (Definition~\ref{D:nice}). For example, $X$ is nice whenever both $X$ and $X^{-1}$ are uniformly bounded or, without uniform boundedness, if $\log X$ has independent increments and $X_{T}$ is in $L^{2}(\P)$ (Proposition~\ref{P:210425}).

It remains to address one final point. There are circumstances when the discounted problem is preferable. This notably happens for Asian options in a model with IID returns, where a suitable choice of numeraire reduces the number of state variables (Ve\v{c}e\v{r} and Xu \cite{vecer.xu.04}). However, there are also situations when solving the discounted problem instead of the direct 
$(S,H,\P)$ minimization is counterproductive. For example, in models with independent asset returns the optimal portfolio depends only on the first and second moments of one-period returns (see, e.g., \cite{li.ng.00}), while the passage to the discounted problem requires full specification of the jump measure for $S$, and not just its first two moments under $\P$. Clearly, in this case the numeraire change introduces a number of unnecessary calculations. 

Motivated by this observation, we use the links in \eqref{eq:3} and \eqref{eq:4} to derive a general solution that deals directly with the undiscounted problem $(S,H,\P)$ in terms of the semimartingale characteristics of $S$ under $\P$. 
We start from the state of knowledge in \cite{cerny.kallsen.07} and \cite{czichowsky.schweizer.13}. We reinterpret those findings suitably (Theorem~\ref{T:main X=1}) and apply them after numeraire change (Corollary~\ref{C:mainhP}). The key novelty of the paper is that we are then able to translate these results back to the original setting without numeraire change (Theorem~\ref{T:main}). On the way, we introduce analogues of the opportunity process (Subsection~\ref{SS:L}) and the mean value process \eqref{eq:V} from \cite{cerny.kallsen.07} in the undiscounted setting. This yields explicit formulae for the optimal trading strategies in terms of semimartingale characteristics of the underlying processes in Theorem \ref{T:explicit}. 

The literature on quadratic hedging in continuous time is largely limited to the formulation \eqref{eq:1}, cf. Schweizer \cite{schweizer.94,schweizer.96}; Bertsimas, Kogan, and Lo~\cite{bertsimas.al.01}; \cite{cerny.kallsen.07}, and \cite{czichowsky.schweizer.13}. Notable exceptions are Zhou and Yin \cite{zhou.yin.03}, Lim~\cite{lim.04,lim.05}, and Yao, Li, and Chen~\cite{yao.al.14}, for example. Notwithstanding these important contributions, Theorems~\ref{T:main} and \ref{T:explicit} and the streamlined efficient frontier formula \eqref{eq:efintro} advance our understanding of the efficient frontier formation in the most general case where a risk-free asset may not be present. Table~\ref{T} below connects existing literature to the three quantities identified in the efficient frontier equation \eqref{eq:efintro}. It is evident that no clear pattern emerges from previous work.

\begin{table}[H] 
\begin{center}
{\footnotesize
\begin{tabular}{lccccc}
\hline\\[-2ex]
 & $L$ & $V(1)$ & $V(H)$ & $\varepsilon ^{2}(1)$ & $\varepsilon ^{2}(H)$\\[1ex] \hline\\[-2ex]
Li and Ng \cite{li.ng.00} & $\tau\ \star$ & $\displaystyle\frac{\mu}{\tau}\ \star$ & --- & $\displaystyle 1- 4\nu^2-\frac{\mu^2}{\tau}\ \star$ & --- \\[1.5ex]  
Bertsimas et al. \cite{bertsimas.al.01}& $a(\cdot,P,Z)$ & from $V(H)$ & $b(\cdot,P,Z)$ & 0 & $c(\cdot,P,Z)$  \\[1ex]  
Zhou and Yin \cite{zhou.yin.03} & $P(\cdot ,i)$ & $H(\cdot ,i)$ & --- & $\theta\ \star$& --- \\[1ex]
Lim \cite{lim.04} & $p$ & $g$ & --- & $M\ \star$ & --- \\[1ex] 
Lim \cite{lim.05} & $p$ & from $V(H)$ & $h$ & from $\eps(H)$ & $\ \star$\tablefootnote{The reference implicitly yields an expression for $\eps^2_0(H)$ but does not give it a label.} \\[1ex] 
Yao et al. \cite{yao.al.14} & $\displaystyle\frac{p}{2}$ & $\displaystyle-\frac{g}{p}$ & --- & $\displaystyle c -\frac{g^2}{2p}$  & --- \\[2ex]
\hline
\end{tabular}
}
\end{center}\medskip
\caption{Notation in pre-existing literature for the key processes $L$, $V(1)$, $V(H)$, $\eps^2(1)$, and $\eps^2(H)$ related to the efficient frontier. Entries marked with $\star$ do not denote a process but only a single value corresponding to $t=0$. Papers ordered chronologically by publication date.\label{T}}
\end{table}

The paper is not the first to treat numeraire change in quadratic hedging. Specific numeraire changes appear in Gourieroux, Laurent, and Pham~\cite{gourieroux.al.98}, Arai~\cite{arai.05}, and Kallsen, Muhle-Karbe, and Vierthauer \cite{kallsen.al.14}; this paper, however, offers the first systematic study of general numeraire changes in the quadratic hedging context. 

\section{Semimartingale characteristics and notation.}\label{S:notation}
A glossary of mathematical notation and terminology is available in Appendix~\ref{App:E}. Superscripts refer generally to coordinates of a vector or vector-valued process rather than powers. The few exceptions should be obvious from the context. If $X$ is a semimartingale, $L(X)$ denotes the set of $X$-integrable predictable processes in the sense of  \cite[III.6.17]{js.03}.

In the subsequent sections, optimal hedging strategies are expressed in terms of semimartingale characteristics.
\begin{definition}[Semimartingale characteristics]\label{d:bcfa}
Let $X$ be an $\R^d$-valued semimartingale with characteristics $(B,C,\nu)$ relative to some truncation function $h:\R^d\to\R^d$.

By \cite[II.2.9]{js.03}, there exists some increasing, predictable, integrable process $A$, some predictable $\R^{d\times d}$-valued process $c$ whose values are non-negative, symmetric matrices, and some transition kernel $F$ from $(\Omega\times\R_+,\P)$ into $(\R^d,\mathscr{B}^d)$ such that
\begin{equation*}
B_t=b\sint A_t, \ \
C_t=c\sint A_t, \ \ 
\nu([0,t]\times G)=F(G)\sint A_t\ \ \mbox{ for }t\in[0,T],\, G\in\mathscr{B}^d.\label{e:bcfa}
\end{equation*}
Here $\mathscr{B}^d$ is the Borel $\sigma$-algebra over $\R^d$. We call $(b,c,F,A)$ {\em differential characteristics} of $X$.\qed
\end{definition}
One should observe that the differential characteristics are not unique: e.g.\ 
$(2b,2c,2F,\frac{1}{2}A)$ yields another version.
Especially for $A_t=t$, one can interpret $b_t$ or rather $b_t+\int(x-h(x))F_t(\d x)$
as a drift rate, $c_t$ as a diffusion coefficient, and $F_t$ as a local jump measure.
The differential characteristics are typically derived from other ``local'' representations of the process, e.g.,
in terms of a stochastic differential equation.

From now on, we choose the same fixed process $A$ for all the (finitely many) semimartingales in this paper.
The results do not depend on its particular choice. In concrete models, $A$ is often taken to be  $A_t=t$ (e.g., for L\'evy processes, diffusions, It\^o semimartingales) or $A_t=[t]=\max\{n\in\N:n\leq t\}$ for discrete-time processes. Since almost all semimartingales of interest in this paper are actually special semimartingales, we use from now on the otherwise forbidden `truncation function'
$h(x)=x,$ which simplifies a number of expressions considerably.

By $\langle X,Y\rangle$ we denote the $\P$-compensator of $[X,Y]$ provided that $X,Y$ are semimartingales such that $[X,Y]$ is $\P$-special (cf.~\cite[p.~37]{jacod.79}). If $X$ and $Y$ are vector-valued, then $[X,Y]$ and $\langle X,Y\rangle$ are to be understood as matrix-valued processes with components $[X^i,Y^j]$ and $\langle X^i,Y^j\rangle$, respectively. If both the integrator $X$ and the integrand $\theta\in L(X)$ are $\R^d$-valued, we define the stochastic integral for vector-valued processes as in  \cite[III.6.17]{js.03}. However, for compatibility with matrix notation, we consider integrators $X$ as column vectors in $\R^d$ and integrands $\theta$ as a row vectors. For $\theta\in L(X)$ and $\psi\in L(Y)$, this yields $[\theta\sint X,\psi \sint Y]=\int \theta \d[X,Y] \psi^\top$ and accordingly $ \langle \theta\sint X,\psi\sint Y\rangle=\int \theta \d\langle X,Y\rangle \psi^\top$. If $\P^\star$ denotes another probability measure, we write $\langle X,Y\rangle^{\P^\star}$ for the $\P^\star$-compensator of $[X,Y]$.

In the whole paper, we write $M^X$ for the local martingale part and $B^X$ for the predictable part of finite variation in the canonical decomposition
$X=X_0+M^X+B^X$ 
of a special semimartingale $X$.
If $\P^\star$ denotes another probability measure, we write accordingly
$X=X_0+M^{X\star}+B^{X\star}$
for the $\P^\star$-canonical decomposition of $X$.

From now on we use the notation $(b^X,c^X,F^X,A)$ to denote differential characteristics of a special semimartingale $X$. If $[X,X]$ is special (i.e., $X$ is locally square-integrable), 
\begin{equation*}
\tilde c^X=c^X+\int xx^\top F(\d x)=b^{[X,X]}
\end{equation*}
stands for the modified second characteristic of $X$. In such case, by \cite[I.4.52]{js.03},  one has 
\[
\langle X,X\rangle = \tilde c^X\sint A.
\]
If they refer to some probability measure $\P^\star$ rather than $\P$, we write instead $(b^{X\star},c^{X\star},F^{X\star},A)$ and
$\tilde c^{X\star}$, respectively.
We denote the joint characteristics of two special vector-valued semimartingales $X,Y$, i.e., the characteristics of $S=(X,Y)$ as
$$\left(b^S,c^S,F^S,A\right) = \left(b^{X,Y},c^{X,Y},F^{X,Y},A\right)= 
\left(\genfrac(){0pt}{0}{b^X}{b^Y},
\left(\begin{array}{cc}
c^{X} & c^{XY} \\ c^{YX} & c^{Y} 
\end{array}\right), F^{X,Y},A\right)$$
and
$$ \tilde c^{X,Y}=\left(\begin{array}{cc}
\tilde c^{X} & \tilde c^{XY} \\ \tilde c^{YX} & \tilde c^{Y} 
\end{array}\right).
$$

In the whole paper, we write $c^{-1}$ for the Moore--Penrose pseudoinverse of a matrix or 
matrix-valued process $c$, which is a particular matrix satisfying $cc^{-1}c=c$ (see Appendix \ref{App:A}).
From the construction, it follows that the mapping $c\mapsto c^{-1}$ is Borel-measurable from $\R^{n\times n}$ to $\R^{n\times n}$ with the Euclidean norm. Moreover, $c^{-1}$ is non-negative and symmetric if this holds for $c$.

\section{Symmetric definition of admissibility.}\label{S:adm}

\subsection{Self-financing strategies.}\label{SS:sf}
Hereafter $S=(X,Y)$ is a  semimartingale in $\R^d$, $d\geq 2$, with $\R$-valued $X$ and 
$\R^{d-1}$-valued $Y$. A trading strategy $\vt\in L(S)$ is accordingly partitioned as $\vt=(\chi,\eta)$. Recall that the integrator $S$ is treated as a column vector and the integrand $\vt\in L(S)$ as a row vector, hence in particular $\vt S=\sum_{i=1}^d\vt^i S^i$. With Harrison and Kreps \cite{harrison.kreps.79}, we use the following numeraire-invariant definition of self-financing trading strategies. Observe that none of the assets in $S$ plays a special role, hence the definition treats all assets \emph{symetrically}.
\begin{definition}[Self-financing strategy]\label{D:SF}
A trading strategy $\vt=(\vt_t)_{0\leq t\leq T} $ in $L(S)$ is called \emph{self-financing} if it satisfies $\vt S=\vt _{0}S_{0}+\vt \sint S $ on $[0,T]$.\qed
\end{definition}
We have the following equivalent characterization of the self-financing condition. 
\begin{proposition}[Equivalent characterizations of self-financing property]\label{P:sf}
The following assertions are equivalent.
\begin{enumerate}[(1)]
\item\label{P:sf.i} $\vt \in L(S)$ and $\vt _{0}S_{0}+\vt \sint S=\vt S$, which implies $(\vt S)_{-}=\vt S_{-}$.

\item\label{P:sf.ii} $\vt \in L(S),$ $\vt _{0}S_{0}+\vt \sint S_{-}=\vt S_{-}$.

\item\label{P:sf.iii} for arbitrary scalar-valued semimartingale $Z$ one has $\vt \in
L(SZ)$ and 
\begin{equation*}
(\vt _{0}S_{0})Z_{0}+\vt \sint (SZ)=(\vt S)Z.
\end{equation*}
\end{enumerate}
\end{proposition}

\bpf
\ref{P:sf.i} 
On taking left limits one obtains
$ \vt _{0}S_{0}+\vt \sint S_- = (\vt S)_-$.
The definition of stochastic integral yields $\vt _{0}S_{0}+\vt \sint S_{-}+\vt \Delta S = \vt S$, 
and on rearranging  $\vt _{0}S_{0}+\vt \sint S_{-} = \vt S_-$. Combined together, one obtains $(\vt S)_{-}=\vt S_{-}$.

\ref{P:sf.ii}\,$\Rightarrow $\,\ref{P:sf.i} follows by adding $\vt\Delta S$ to both sides of \ref{P:sf.ii}. 

\ref{P:sf.i}\,$\Rightarrow$\,\ref{P:sf.ii} follows trivially by taking left limits.
To show \ref{P:sf.iii}\,$\Rightarrow $\,\ref{P:sf.i}, take $Z=1$. 

\ref{P:sf.i}\,$\Rightarrow$\,\ref{P:sf.iii} Integration by parts yields
\begin{align}
(\vt S)Z &{}=(\vt S)_{0}Z_{0}+(\vt S)_{-}\sint Z+Z_{-}\sint (\vt S)+[\vt S,Z]  \notag \\
&{}=(\vt S)_{0}Z_{0}+(\vt S)_{-}\sint Z+Z_{-}\sint (\vt \sint S)+[\vt \sint S,Z].  \label{sf5}
\end{align}
From \ref{P:sf.i} one has $(\vt S)_{-}\sint Z=(\vt S_{-})\sint Z$. Goll and Kallsen~\cite[Proposition 5.1]{goll.kallsen.03} now yields
$\vt \in L(S_{-}\sint Z)$. Similarly, by \cite[Proposition 5.1 and 5.2]{goll.kallsen.03} $\vt \in L(Z_{-}\sint S)$ and 
$\vt \in L([S,Z])$. Integrating by parts, $SZ=S_{0}Z_{0}+S_{-}\sint Z+Z_{-}\sint S$, one obtains $\vt \in L(SZ)$. After simplifications, \eqref{sf5} yields
$(\vt S)Z=(\vt S)_{0}Z_{0}+\vt \sint (SZ)$.
\ep
\begin{definition}[Self-financing strategy on a subinterval]\label{D:SF sub}
For a $[0,T]$-valued stopping time $\tau$, a process $\vt=(\vt_t)_{0\leq t\leq T}$ is a \emph{self-financing trading strategy on} $\auf\tau,T\zu$ 
if $\one_{\zu\tau,T\zu}\vt\in L(S)$, $\vt_\tau$ is $\F_\tau$-measurable, and $\vt S=\vt_\tau S_\tau+(\one_{\zu\tau,T\zu}\vt)\sint S$ on $\auf\tau,T\zu$. \qed
\end{definition}

\begin{remark}[Self-financing strategy on {$\auf\tau,T\zu$} vs. {$[0,T]$}]
Observe that, for any $[0,T]$-valued stopping time $\tau$, a self-financing strategy on $[0,T]$ is also self-financing on $\auf\tau,T\zu$. Moreover, Proposition \ref{P:sf} holds on $\auf\tau,T\zu$. Note, however, that in Definition~\ref{D:SF sub} $\vt_\tau$ is only $\F_\tau$-measurable, while $\vt_\tau$ of a self-financing strategy on $[0,T]$ must be $\F_{\tau-}$-measurable. \qed
\end{remark}

\subsection{Admissibility and \texorpdfstring{$L^2(\P)$}{L2(P)}--closedness.}

Following \cite{biagini.cerny.20}, admissible strategies are obtained as an appropriate closure of `tame strategies' whose wealth is highly integrable. Unlike \cite{cerny.kallsen.07}, we do not insist tame strategies are of buy-and-hold type; this yields better properties when the price process is enlarged by a self-financing strategy (see Remark~\ref{R:210425} below). 

\begin{definition}[Tame trading strategies]\label{D:L2P}
\hspace*{-1em} We say that $\vt=(\vt_t)_{0\leq t\leq T}$ is a \emph{tame trading strategy}, writing $\vt \in \tame (S,\P)$, if
\begin{enumerate}[(i)]
\item\label{D:L2P.1} $\vt $ is self-financing, i.e., it satisfies $\vt S=\vt _{0}S_{0}+\vt \sint S $ on [0,T] with $\vt\in L(S)$;
\item\label{D:L2P.2} $\vt S$ is an $L^{2}(\P)$-semimartingale in the sense of Delbaen and Schachermayer \cite{delbaen.schachermayer.96.aihp}, i.e., 
$$\sup \{\E[(\vt _{\sigma }S_{\sigma})^{2}]~:~\sigma \text{ is a $[0, T]$-valued stopping time}\}<\infty .$$
\end{enumerate}
For $v\in L^{0}(\F_{0},\P)$, we let $\tame_{v}(S,\P)=\{\vt \in \tame (S,\P)~:~\vt _{0}S_{0}=v\}$.\qed
\end{definition}

\begin{remark}[Invariance of tame strategies to self-financed market extension]\label{R:210425}
Observe that one can take any self-financed wealth process and add it as an extra component of the price process $S$ without altering the set of tame terminal wealths. That is,  for \emph{any} self-financing $\vp \in L(S)$ we automatically have that
\begin{equation*}
\left\{ \vt _{T}S_{T}~:~\vt \in \tame _{0}(S,\P)\right\} 
=\left\{\theta _{T}(S_{T},\vp _{T}S_{T})~:~\theta \in \tame _{0}((S,\vp S),\P)\right\} .\tag*{\qed}
\end{equation*}
\end{remark}

\begin{remark}[Intermediate wealth]
Note that when there is no risk-free asset, the quantity $\vt _{\tau }S_{\tau}$ appearing in Definition~\ref{D:L2P}\ref{D:L2P.2}, cannot be interpreted as the terminal value of a trading strategy liquidated at time $\tau$. That is, there will typically be no self-financing strategy that turns the wealth $\vt _{\tau }S_{\tau}$ at time $\tau$ into the wealth $\vt _{\tau }S_{\tau}$ at time $T$.\qed
\end{remark}
\begin{remark}[Tame strategies and non-existence of solutions]
In finite discrete time with a constant risk-free asset, the intermediate wealth of an optimal strategy is not necessarily in $L^2(\P)$, hence the set of tame terminal wealths is generally not closed in $L^2(\P)$; see \citep[Example~8.6]{cerny.kallsen.09} and related continuous-time results of  Stricker and co-authors~\cite{choulli.al.98,delbaen.al.97,monat.stricker.95}. In a general setting with constant risk-free asset, \v{C}ern\'y and Kallsen \cite{cerny.kallsen.07} show how to obtain a closed set of terminal wealths by passing to a wider class of `admissible integrands.' This idea is extended below to the case without a risk-free asset. \qed
\end{remark}

\begin{definition}[Admissible trading strategies]\label{D:adm}
We say that $\vt$ is an \emph{admissible} trading strategy, writing $\vt \in \adm(S,\P)$, if 
$\vt$ is self-financing in the sense of Definition~\ref{D:SF} and
there is an approximating sequence of tame trading strategies $\vt ^{(n)}\in \tame(S,\P)$ such that
\begin{enumerate}[(i)]
\item $\vt _{\sigma}^{(n)}S_{\sigma}\overset{\P}{\rightarrow }\vt _{\sigma}S_{\sigma}$ for all $[0,T]$-valued stopping times $\sigma$;
\item $\vt _{T}^{(n)}S_{T}\overset{L^{2}(\P)}{\rightarrow }\vt_{T}S_{T}$.
\end{enumerate}
For $v\in L^{0}(\F_{0},\P)$, we set $\adm_{v}(S,\P)=\{\vt \in \adm (S,\P)~:~\vt _{0}S_{0}=v\}$.\qed
\end{definition}

We next give sufficient conditions for the admissible terminal wealths to be closed in $L^{2}(\P)$. 

\begin{definition}[Deflator]\label{D:deflator}
We call process $Z=(Z_t)_{0\leq t\leq T}$ a \emph{deflator} for $S$ under $\P$, writing $Z\in \cM_{2}(S,\P)$, if
\begin{enumerate}[(i)]
\item\label{D:deflator.i} $Z_{T}\in L^{2}(\P)$;
\item\label{D:deflator.ii} $Z(\vt S)$ is a $\P$-martingale for all $\vt \in \tame _{0}(S,\P)$;
\item\label{D:deflator.iii} $Z(\vp S)$ is a $\P$-martingale for some self-financing strategy $\vp$ such that $\vp S>0$.\qed
\end{enumerate}
\end{definition}

\begin{lemma}[Martingale property of deflated admissible wealth]\label{lem: (theta S)Z}
For all admissible strategies $\vt \in \adm(S,\P)$ and all deflators $Z\in\cM_2(S,\P)$, the process $Z(\vt S)$ is a $\P$-martingale.
\end{lemma}
\bpf
Consider $\vt \in \adm(S,\P)$ and the corresponding sequence of approximating simple strategies $\vt ^{(n)}\in \tame(S,\P)$. Since $\vt _{T}^{(n)}S_{T}\overset{L^{2}(\P)}{\rightarrow }\vt _{T}S_{T}$, the H\"{o}lder inequality yields  
\begin{equation}\label{ZSP1}
\E[\vt _{T}^{(n)}S_{T}Z_{T}|\F_{t}]\overset{L^{1}(\P)}\rightarrow \E[\vt _{T}S_{T}Z_{T}|\F_{t}].
\end{equation}
The martingale property of $Z(\vt ^{(n)}S)$ in Definition~\ref{D:deflator}\ref{D:deflator.ii} yields $\E[Z_{T}(\vt _{T}^{(n)}S_{T})|\F_{t}]=Z_{t}(\vt_{t}^{(n)}S_{t})$, while from the definition of admissibility $Z_{t}(\vt _{t}^{(n)}S_{t})\overset{\P}{\rightarrow }Z_{t}(\vt_{t}S_{t})$ and as a consequence of \eqref{ZSP1} we obtain 
$\E[Z_{T}(\vt_{T}S_{T})|\F_{t}]=Z_{t}(\vt _{t}S_{t})$, $\P$-almost surely.
\ep

\begin{theorem}[Existence\hspace{-0.1em} and\hspace{-0.1em} uniqueness\hspace{-0.1em} of\hspace{-0.1em} the\hspace{-0.1em} optimal\hspace{-0.1em} wealth\hspace{-0.1em} process]\label{T:closedness}\!\!\!
If $\cM_2(S,\P)$  contains a strictly positive element, the hedging problem \eqref{eq:1new} allows for an optimizer 
$\vt\in\adm_v(S,\P)$, provided $\adm_v(S,\P)$ is non-empty. Moreover, the value process $\vt S$ of the optimizer is up to indistinguishability unique.
\end{theorem}
\bpf
In steps \ref{pTclosedness.i}--\ref{pTclosedness.vi} below we shall show, for each $v\in \R$, that the set 
$$K_{v}(S,\P)=\{\vt _{T}S_{T}:\vt \in \adm_{v}(S,\P)\},$$ 
if not empty, is a closed affine subspace of $L^{2}(\P)$.
Step \ref{pTclosedness.vii} then argues existence and uniqueness of the optimal wealth process.
\begin{enumerate}[(i),wide = 0pt]
\item\label{pTclosedness.i} Assume $K_{1}(S,\P)$ is non-empty, otherwise there is nothing to prove for $v\neq 0$. Since $\adm(S,\P)$ is a vector space one has $K_{v}(S,\P)=vJ+K_{0}(S,\P)$ for some $J\in K_{1}(S,\P)$ and it is enough to prove that $K_{0}(S,\P)$ is closed. Consider a random variable $G\in L^{2}(\P)$ and a convergent sequence
$$\vt_{T}^{(n)}S_{T}\overset{L^{2}(\P)}{\longrightarrow }G,\qquad \vt^{(n)}\in \adm_{0}(S,\P).$$
\item\label{pTclosedness.ii} By assumption there is $Z>0$ in $\cM_2(S,\P)$. By Lemma \ref{lem: (theta S)Z}, $Z(\vt ^{(n)}S)$ is a $\P$-martingale for each $n\in\N$. Denoting $N=\vp S>0$, we have by assumption that $NZ>0$, too, is a $\P$-martingale, hence $N_->0$. Without loss of generality, we may assume $X\equiv S^{1}=N$ because adding a self-financed price process $N=\vp S$ to $S$ does not affect the set of terminal wealths attained by tame strategies; see Remark~\ref{R:210425}. On defining a new measure 
$\tP\sim \P$ by letting
\begin{equation*}
\frac{\d\tP}{\d\P}=\frac{X_{T}Z_{T}}{X_{0}Z_{0}},
\end{equation*}
the martingale property of $XZ$ yields the density process of the form 
\begin{equation}
\frac{\d(\tP|_{\F_{t}})}{\d(\P|_{\F_{t}})}=\E\left( \frac{\d\tP}{\d\P}\,\bigg|\,\F_{t}\right) =\frac{X_{t}Z_{t}}{X_{0}Z_{0}}.
\label{eq: Phat_density}
\end{equation}
\item\label{pTclosedness.iii} Let $\hS=\frac{S}{X}=(1,\hY)$. Process $\vt ^{(n)}\hS$ is a $\tP$-martingale by \cite[III.3.8]{js.03} and \eqref{eq: Phat_density} .
\item\label{pTclosedness.iv} H\"{o}lder's inequality yields $\vt_T ^{(n)}S_{T}Z_{T}\overset{L^{1}(\P)}{\longrightarrow }GZ_{T}$, equivalent to the $L^{1}(\tP)$-convergence of $(\vt ^{(n)}\hS)_{T}$ to the random variable $\hG =\frac{G}{X_{T}}$.
\item\label{pTclosedness.v} Then, $\vt ^{(n)}\in \adm _{0}(S,\P)$ and Proposition~\ref{P:sf} yield $\vt ^{(n)}\hS=\vt ^{(n)}\sint \hS$. In view of \ref{pTclosedness.iii} and \ref{pTclosedness.iv}, a classical martingale compactness result of Yor~\cite[Corolaire 2.5.2]{yor.78} yields the existence of $\eta \in L(\hY,\tP)$ such that $\eta \sint \hY$ is a $\tP$-martingale and $\eta \sint \hY_{T}=\hG $. We now let 
\begin{equation*}
\chi =\eta \sint \hY-\eta \hY=(\eta \sint \hY)_{-}-\eta \hY_{-}.  
\end{equation*}
By construction $\chi $ is predictable and therefore it belongs to $L(\hS^{1})= L(1)$. Consequently $\vt =(\chi ,\eta )$ is in $L(\hS) $ and by construction $\vt $ is self-financing; 
\begin{equation*}
\vt \sint \hS=\eta \sint \hY=\chi +\eta \hY= \vt \hS.  
\end{equation*}
\item\label{pTclosedness.vi} The $L^{1}(\tP )$-convergence of $\vt _{T}^{(k)}\hS_{T}$ to $\vt _{T}\hS_{T}= \hG $ established in \ref{pTclosedness.iv} yields 
\begin{equation*}
\vt _{t}^{(k)}\hS_{t}
=\E^{\tP}[(\vt _{T}^{(n)}\hS_{T})|\F_{t})\overset{L^{1}(\tP)}{\rightarrow }\E^{\tP}((\vt _{T}\hS_{T})|\F_{t})
=\vt _{t}\hS_{t}.
\end{equation*}
In view of $\P\sim \tP$ and $X>0$ this yields $\vt_{t}^{(k)}S_{t}=(\vt_{t}^{(k)}\hS_{t})X_t\overset{\P}{\rightarrow }(\vt _{t}\hS_{t})X_t=\vt _{t}S_{t}$, which completes
the proof of $L^2(\P)$-closedness.
\item\label{pTclosedness.vii} $L^2(\P)$-closedness of $K_v(S,\P)$ and strict convexity in \eqref{eq:1new} yield the existence of an optimizer together with $\vt_T S_T=\tilde\vt_T S_T$
almost surely for any two optimizers $\vt,\tilde\vt$. The martingale property in Lemma~\ref{lem: (theta S)Z}  yields
\[Z_t\vt_t S_t=\E[Z_T\vt_T S_T|\F_t]=\E[Z_T\tilde\vt_T S_T|\F_t]
=Z_t\tilde\vt_t S_t\]
and hence $\vt_t S_t=\tilde\vt_t S_t$ almost surely because $Z$ is positive.\qed
\end{enumerate}

\begin{remark}[On the assumption of Theorem~\ref{T:closedness}]\label{R:deflators}
\begin{enumerate*}[(1)]
\item\label{R:deflators1} \ Since $Z(\vp S)$ in Definition~\ref{D:deflator} is a strictly positive $\P$-martingale, neither the processes $Z$ and $\vp S$ nor their left limits $Z_-$ and $(\vp S)_-$ are allowed to hit zero. One could weaken these assumptions by considering $Z$ and $\vp S$ in the form of stochastic exponentials that are restarted after hitting zero, \emph{\`a la} Choulli, Krawczyk, and Stricker~\cite{choulli.al.98}. Proposition 6.1 and Theorem 6.2 of \cite{czichowsky.schweizer.13} then offer a way of generalizing Theorem~\ref{T:closedness} to such relaxed setting; see \cite{cerny.czichowsky.25} for the details.\smallskip\\ 
\item\label{R:deflators2} \ For the assumption of Theorem~\ref{T:closedness} to hold, it is sufficient that there exists a strictly positive process $Z$ such that $ZS$ is a $\P$-$\sigma$-martingale and $\sup_{t\in \lbrack 0,T]}|Z_{t}|\in L^{2}(\P)$. Up to the $L^2(\P)$-condition, this holds if $X>0$ and $\hS=S/X$ satisfies the no free lunch with vanishing risk (NFLVR) condition, which is implied by the existence of an equivalent local martingale measure for $\hS$; see Delbaen and Schachermayer~\cite{delbaen.schachermayer.98}.
In that sense, the assumption in Theorem~\ref{T:closedness} is related to the absence of arbitrage.\qed
\end{enumerate*}
\end{remark}

\subsection{Opportunity process.}\label{SS:L}
The next definition describes a straightforward extension of tameness and admissibility to a subinterval. It is needed to establish the key concept of the opportunity process below.

\begin{definition}[Tame and admissible trading strategies on a subinterval]\label{D:tame adm sub}
\ \ For a $[0,T]$-valued stopping time $\tau$, we say that $\vt=(\vt_t)_{0\leq t\leq T}$ is a \emph{tame trading strategy on} $\auf \tau,T\zu$, 
writing $\vt \in \tame^\tau (S,\P)$, if
\begin{enumerate}[(i)]
\item $\vt $ is self-financing on $\auf \tau,T\zu$ in the sense of Definition~\ref{D:SF sub}; 

\item $\esssup \{\E[(\vt _{\sigma }S_{\sigma})^{2}|\F_\tau]~:~ \sigma \text{ is a $[\tau,T]$-valued stopping time}\}<\infty$, where we use the generalized conditional expectation as in Jacod and Shiryaev \cite[I.1.1]{js.03}.
\end{enumerate}

We say that $\vt$ is an admissible trading strategy on $\auf \tau,T\zu$, writing $\vt \in \adm^{\,\tau}(S,\P)$, if 
$\vt$ is self-financing on $\auf \tau,T\zu$ in the sense of Definition~\ref{D:SF sub} and there is an approximating sequence of tame trading strategies $\vt ^{(n)}\in \tame^\tau(S,\P)$ such that
\begin{enumerate}[(i),resume]
\item $\vt _{\sigma}^{(n)}S_{\sigma}\overset{\P}{\rightarrow }\vt _{\sigma}S_{\sigma}$ for all $[\tau,T]$-valued stopping times $\sigma$;
\item $\vt _{T}^{(n)}S_{T}\rightarrow\vt_{T}S_{T}$ in $L^{2}(\P[\,\cdot\,|\F_{\tau}])$.
\end{enumerate}
Here, a sequence $X_n$ of random variables converges to a random variable $X$ in $L^{2}(\P[\,\cdot\,|\F_{\tau}])$ if the sequence of random variables $(\E[|X_n-X|^2|\F_{\tau}])_{n=1}^\infty$ converges to $0$, $\P$-almost surely.

Finally, for $v\in L^{0}(\F_{\tau},\P)$, we let
\begin{equation*}
\tame^\tau_{v}(S,\P)=\{\vt \in \tame^\tau (S,\P)~:~\vt _{\tau}S_{\tau}=v\};\qquad 
\adm^{\tau}_{v}(S,\P)=\{\vt \in \adm^{\tau} (S,\P)~:~\vt _{\tau}S_{\tau}=v\}.\tag*{\qed}
\end{equation*}
\end{definition}

We shall now extend the notion of the \emph{opportunity process} from \cite{cerny.kallsen.07} to the setting without a risk-free asset. This extension is non-trivial: Corollary~3.4 in \cite{cerny.kallsen.07} interprets the opportunity process as the conditional squared hedging error from approximating the constant payoff 1 by portfolios that cost 0, which is then linked to the Sharpe ratio and therefore to investment opportunities. What is needed instead is the squared hedging error in the approximation of the constant payoff 0 by portfolios that cost 1 (called the fully invested portfolios). The two quantities happen to coincide when there is a risk-free asset with constant value 1 but in the general setting studied here their roles are quite different. The latter no longer has a direct link to the Sharpe ratio of zero-cost portfolios; instead, it represents the minimal conditional second moment among fully invested portfolios. Thus the terminology ``opportunity process,'' which we shall maintain, is a misnomer in the context of this paper.    

\begin{definition}[Opportunity process] The \emph{opportunity process} $L=(L_t)_{0\leq t\leq T}$ is given by
\begin{equation} \label{D:op}
L_t=\underset{\vt\in\adm^{t}_{1}(S,\P)}{\essinf}\E\left[(\vt_TS_T)^2|\F_t\right]=\underset{\vt\in\adm^{t}_{1}(S,\P)}{\essinf}\E\left[(1+(\one_{\zu t,T\zu}\vt)\sint S_T)^2\big|\F_t\right],\quad 0\leq t\leq T.
\end{equation}
Thus, $L_t$ measures the smallest conditional second moment among fully invested self-financing portfolios on the subinterval $[t,T]$.\qed
\end{definition}
\subsection{Nice numeraire.}\label{SS:5.1}
For the proof of our explicit representation of optimal strategies, the existence of sufficiently well-behaved numeraires plays an important role.

\begin{definition}[Nice numeraire]\label{D:nice}
We say that $N$ is a \emph{nice numeraire} under $\P$ if $N>0$, $N_->0$, and there are constants $\underline{\delta}$, $\overline{\delta}$ such that, for all $t\in \lbrack 0,T]$, we have
\begin{equation}\label{eq:vnice}
0<\underline{\delta} \leq \frac{\E[N_{T}^{2}|\F_{t}]}{N_{t}^{2}}\leq \overline{\delta}<\infty .
\end{equation}
We say that $S$ \emph{admits a nice numeraire} if there is a self-financing strategy $\vp$ for $S$ such that $\vp S$ is a nice numeraire.
\qed
\end{definition}
Recall the notion of the deflator in Definition~\ref{D:deflator}.
\begin{proposition}[Consequences of nice numeraire existence]\label{P:ThetaEqual}
Suppose $X\equiv S^1$ is a nice numeraire. Then,
\begin{enumerate}[(1)]
\item\label{ThetaEqual.1} $\tame(S,\P)=\tame(\hS,\hP)$ and $\tame^\tau(S,\P) = \tame^\tau(\hS,\hP)$ for any $[0,T]$-valued stopping time $\tau$;
\item\label{ThetaEqual.1bis} $S$ is locally square-integrable under $\P$ if and only if $\hS$ is locally square-integrable under $\hP$;
\item\label{ThetaEqual.2} $X$ is the wealth of a tame strategy for $S$ under $\P$;
\item\label{ThetaEqual.3} $ZX$ is a $\P$-martingale for any $Z\in\cM_2(S,\P)$;
\item\label{ThetaEqual.4} $ZX/\hZ\in \cM_2(\hS,\hP)$ if and only if $Z\in\cM_2(S,\P)$, where $\hZ_t=\E[X_T^2|\F_t]$. 
\end{enumerate}
\end{proposition}
\bpf
\ref{ThetaEqual.1} We only give the proof for $\tame(S,\P) = \tame(\hS,\hP)$; that for $\tame^\tau(S,\P) = \tame^\tau(\hS,\hP)$ for any $[0,T]$-valued stopping time $\tau$ follows by a straightforward modification of the arguments.
\begin{enumerate}[(i),wide = 0pt]
\item\label{pPThetaBar.i} Since $X$ is nice, by Proposition 3.3 in \cite{choulli.al.98} one has $\E[(\frac{X_{T}}{X_{\tau }})^{2}|\F_{\tau })\leq \overline{\delta}$ for all stopping times $\tau \leq T.$ Additionally,
\begin{equation}\label{eq: reverseHolderX^-1}
\E^{\hP}[( X_{T}^{-1}/X_{\tau }^{-1}) ^{2}|\F_{\tau}]=\E[X_{\tau }^{2}|\F_{\tau }]/\E[X_{T}^{2}|\F_{\tau}]
=1/\E[\left( X_{T}/X_{\tau }\right) ^{2}|\F_{\tau }].
\end{equation}
In view of \eqref{eq: reverseHolderX^-1}, Proposition~3.3 of \cite{choulli.al.98} applied to $X^{-1}$ under $\hP$ yields that for all $[0,T]$-valued stopping times $\tau$ one has $\E[\left( X_{T}/X_{\tau }\right)^{2}|\F_{\tau }]\geq \underline{\delta}$. 

\item\label{pPThetaBar.ii} Consider an $L^{2}(\P)$-semimartingale $W$. The identity 
\begin{align*}
\E[X_{T}^{2}]\E^{\hP}[(W/X)_{\tau }^{2}] &{}=\E[X_{T}^{2}(W/X)_{\tau }^{2}]
{}=\E[W_{\tau }^{2}\E[\left( X_{T}/X_{\tau }\right) ^{2}|\F_{\tau }]]
\end{align*}
and \ref{pPThetaBar.i} yield 
$
\underline{\delta} \E[W_{\tau }^{2}]\leq \E[X_{T}^{2}]\E^{\hP}[(W/X)_{\tau}^{2}]\leq \overline{\delta}\E[W_{\tau }^{2}],
$
for all $[0,T]$-valued stopping times $\tau$. Hence, $W$ is an $L^{2}(\P)$-semimartingale if and only if $\hW =W/X$ is an $L^{2}(\hP)$-semimartingale. This shows that the tame strategies for $(S,\P)$ and $(\hS,\hP)$ coincide.
\end{enumerate}

\ref{ThetaEqual.1bis} This follows from step \ref{pPThetaBar.ii} above by localization since local $L^{2}(\P)$-semimartingales coincide with locally square-integrable semimartingales under $\P$ (and likewise for $\hP$ in place of $\P$) by \cite[Lemma~A.2]{cerny.kallsen.07}.
\smallskip

\ref{ThetaEqual.2}--\ref{ThetaEqual.3} Let $\vp = (1,0,\ldots,0)$. Since $\vp \hS = 1$, we have $\vp \in\tame(\hS,\hP)$. The claims now follow from \ref{ThetaEqual.1} together with Definition~\ref{D:nice}.
\smallskip

\ref{ThetaEqual.4} By \cite[Proposition~III.3.8]{js.03}, $Z(\vt S)$ is a $\P$-martingale if and only if $Z(\vt S)/\hZ$ is a $\hP$-martingale. Since $S=X\hS$, the claim follows.
\ep

We shall now state an easily verifiable sufficient condition for $X$ to be nice.

\begin{proposition}[Nice numeraire with independent returns]\label{P:210425}
Suppose that $X>0$, $X_->0$ and that $\Log(X)$ is a semimartingale with independent increments satisfying $\Log(X)_T\in L^2(\P)$. Then $X$ is a nice numeraire.
\end{proposition}
\bpf
By \cite[Theorem II.2.29]{js.03}, one obtains that $\Log(X)$ and $[\Log(X)]$ are special. \v{C}ern\'{y} and Ruf~\cite[Proposition~2.15 and Theorem~4.1]{cerny.ruf.23.spa} and Yor's formula then yield for $W=2\Log(X) + [\Log(X)]$
\begin{align*}
\frac{\E[X^2_T|\F_t]}{X^2_t}&=\E\left[\Exp(\one_{(t,T]}\sint\Log(X))_T^2\right]
=\E[\Exp(\one_{(t,T]}\sint W)_T]=\Exp\big(\one_{(t,T]}\sint B^{W}\big)_T,\qquad t\in [0,T].
\end{align*}
Since $B^W$ is deterministic and of finite variation with $\Delta B^W>-1$, the claim follows.
\ep

The next example shows that the stochastic exponential of an It\^o processes with bounded coefficients is a nice numeraire. These processes arise in applications of linear quadratic control and BSDE techniques to mean--variance portfolio selection.
\begin{example}[It\^o process with bounded coefficients yields a nice numeraire]\leavevmode\newline Suppose that $S$ is an It\^o process satisfying  
$\d S_t=S_t(\mu_t\d t+\sigma_t\d W_t)$ for $t\in[0,T]$, where
$\mu$ and $\sigma$ are bounded predictable processes and $W=(W_t)_{0\leq t\leq T}$ is an $1$-dimensional Brownian motion. For $A_t=t$, the semimartingale characteristics are then given by $b^{\Log(S)}=\mu$, $c^{\Log(S)}=\sigma^2$ and $F^{\Log(S)}\equiv0$. In this setting, we shall verify that $S$ is a nice numeraire.

Indeed, one has 
\begin{align*}
\frac{\E[S^2_T|\F_t]}{S^2_t}&{}=\E\left[\left(\e^{\int_t^T\left(\mu_s-\frac{1}{2}\sigma^2_s\right)\d s+\int_t^T\sigma_s \d W_s}\right)^2\,\bigg|\,\F_t\right]\\[0.25ex]
&{}=\E\left[\e^{\int_t^T\left(2\mu_s+\sigma^2_s\right)\d s-\frac{1}{2}\int_t^T(2\sigma_s)^2\d s+\int_t^T 2\sigma_s \d W_s}\,\Big|\,\F_t\right]
\leq\left\|\e^{\int_t^T\left(2\mu_s+\sigma^2_s\right)\d s}\right\|_{L^\infty(\P)},
\end{align*}
thanks to H\"older's inequality and $\E[\e^{-\frac{1}{2}\int_t^T(2\sigma_s)^2\d s+\int_t^T 2\sigma_s \d W_s}\,|\,\F_t]=1$. 
By Jensen's inequality, one similarly obtains
\begin{align*}
\frac{S^2_t}{\E[S^2_T|\F_t]}&{}
\leq\E\left[\left(\e^{\int_t^T\left(\mu_s-\frac{1}{2}\sigma^2_s\right)\d s+\int_t^T\sigma_s \d W_s}\right)^{-2}\,\bigg|\,\F_t\right]\\[0.25ex]
&{}=\E\left[\e^{\int_t^T\left(-2\mu_s+3\sigma^2_s\right)\d s-\frac{1}{2}\int_t^T(2\sigma_s)^2\d s-\int_t^T 2\sigma_s \d W_s}\,\Big|\,\F_t\right]\leq\left\|\e^{\int_t^T\left(-2\mu_s+3\sigma^2_s\right)\d s}\right\|_{L^\infty(\P)}.
\end{align*}
Since $\mu$ and $\sigma^2$ are bounded, this yields constants $\underline{\delta},\overline{\delta}$ such that \eqref{eq:vnice} holds for all $t\in[0,T]$ with $N$ replaced by $S$.\qed
\end{example}
\subsection{Compatibility of the symmetric formulation with previous studies.}\label{SS:X=1}
Theorem~\ref{T:main X=1} below recasts the results of \cite{cerny.kallsen.07,czichowsky.schweizer.13} in the symmetric form. It is of independent interest to specialists since it provides a new variational characterization of the so-called adjustment process and other important quantities that appear in the classical setting with $S^{1}=X=1$. However, in the context of this paper, it plays an auxiliary role as a key step in proving the main result of the paper, Theorem~\ref{T:main}.

Recall $Y= S^{2:d}$ and that self-financing strategies $\vt $ for $S=(X,Y)$ are partitioned $\vt = (\chi ,\eta )$.
\begin{proposition}[Transposition of admissibility from \cite{cerny.kallsen.07}]\label{P:strategies}
Suppose that $Y$ is an $\mathbb{R}^{d-1}$-valued locally square-integrable semimartingale admitting an equivalent martingale measure with square-integrable density. Denote by $\adm^{\CK}(Y,\P)$ the set of admissible strategies with zero initial wealth in the sense of \cite{cerny.kallsen.07}. Then,
$$
\adm_{0}(S,\P)=\left\{\vt=(\chi,\eta)~:~\eta\in\adm_{0}^{\CK}(Y,\P)\text{ with $\chi=\eta\sint Y_{-}-\eta Y_{-}$}\right\}
$$
and $\vt S=v+\eta\sint Y$ for all $\vt=(\chi,\eta)\in\adm_{v}(S,\P)$.
\end{proposition}
\bpf
We denote the set of all reduced form strategies corresponding to initial wealth $v$ by $\adm_{v}^{2:d}(S,\P)=\{\vt^{2:d}~:~\vt \in \adm_{v}(S,\P)\}$.
\begin{enumerate}[(i),wide=0pt]
\item\label{pPstrategies:i} Inclusion $\tame^{\CK}(Y,\P)\subseteq \tame_{0}^{2:d}(S,\P)$. Consider $\vt ^{2:d}=\vt_{\sigma}^{2:d}1_{\zu\sigma ,\tau \zu}\in \tame ^{\CK}(Y,\P)$. Define now 
$\vt ^{1}=\vt_{\sigma}^{2:d}(S_{\tau }-S_{\sigma -})1_{\zu\sigma ,T\zu}.$ 
Then, $\vt \in \tame (S,\P)$ since $\vt \sint S=\vt S$, $\vt _{0}S_{0}=0$, and $\vt \sint S=\vt ^{2:d}\sint S^{2:d}$ is an $L^{2}(\P)$-semimartingale.

\item\label{pPstrategies:ii} Inclusion $\tame _{0}^{2:d}(S,\P)\subseteq \adm^{\CK}(S^{2:d},\P)$. The process $\vt ^{2:d}\sint S^{2:d}=\vt \sint S$ is an $L^{2}(\P)$-semimartingale. By Lemma \ref{lem: (theta S)Z}, $(\vt^{2:d}\sint S^{2:d})Z$ is a martingale for every 
$Z\in \cM_2(S,\P)$. By \cite[Corollary~2.5]{cerny.kallsen.07} we have $\vt ^{2:d}\in\adm^{\CK}(S^{2:d},\P)$.

\item\label{pPstrategies:iii} From \ref{pPstrategies:i} and \ref{pPstrategies:ii} we know $\tame ^{\CK}(S^{2:d},\P)\subseteq \tame_{0}^{2:d}(S,\P)\subseteq \adm^{\CK}(S^{2:d},\P)$. On taking closures we obtain $\adm^{\CK}(S^{2:d},\P)\subseteq \adm^{2:d}(S,\P)\subseteq \adm^{\CK}(S^{2:d},\P)$.\ep
\end{enumerate}

On the way to establishing Theorem~\ref{T:main}, we start with the case $X=1$. When reading Theorem~\ref{T:main}, observe that $X=1$ is trivially a nice numeraire in the sense of Definition~\ref{D:nice}. Moreover, with $S=(1,Y)$, a positive element of $\cM_2(S,\P)$ is simply the (square integrable) density process of an equivalent local martingale measure for $Y$; see Definition~\ref{D:deflator}. Observe also that the proof of Theorem~\ref{T:main X=1} does not involve a numeraire change. We do not reproduce the lengthy statement of Theorem~\ref{T:main} here to avoid repetition.

\begin{theorem}[Results of \cite{cerny.kallsen.07,czichowsky.schweizer.13} in a symmetric form]\label{T:main X=1}
Theorem~\ref{T:main} holds if $X=1$.
\end{theorem}
\bpf
See Appendix~\ref{S:proof X=1}.
\ep
\subsection{Optimal hedging after numeraire change.}\label{SS:hX=1}
As yet, we do not know how to obtain the solution of the hedging problem \eqref{eq:1new} for general $X$. 
The discussion in the introduction suggests that the optimizer for the price process $S$, payoff $H$, and probability measure $\P$ 
coincides with the one for the ``discounted model'' with the price process $\hS=(1,\hY)$, payoff $\hH=H/X_T$ and probability measure $\hP$  as in \eqref{eq:hP}. Therefore, it makes sense to apply Theorem~\ref{T:main X=1} to this discounted setup.

\begin{corollary}[Quadratic hedging after numeraire change]\label{C:mainhP}
Suppose that $S$ admits a nice numeraire and $\mathcal M_2(S,\mathsf{P})$ contains a strictly positive element.
Then $\hat S$ is a locally square-integrable semimartingale under $\hat P$ and $\mathcal M_2(\hat S,\mathsf{\hat P})$ contains a strictly positive element.
Moreover, Theorem~\ref{T:main} applies to $\hat S, \hat H, \hat P$, i.e., 
\begin{enumerate}[(1)]
\item\label{C:mainhP.1} The opportunity process $\hL$ is the unique bounded semimartingale $\hL=(\hL_t)_{0\leq t\leq T}$ such that
\begin{enumerate}[(a)]
\item\label{C:mainhP.1a} $\hL>0$, $\hL_->0$, and $\hL_T=1$; 
\item $\frac{\hL}{\Exp(B^{\Log(\hL),\hP})}>0$ is a martingale on $[0,T]$;
\item $[\hS,\hS]$ is $\hP^\star$-special for the measure $\hP^\star\sim \hP$ defined by
$$\frac{\d\hP^\star}{\d\hP}=\frac{\hL_T}{\E^{\hP}[\hL_0]\Exp(B^{\Log(\hL),\hP})_T}>0,$$
which implies
\begin{align}
b^{\hS,\hP^\star}  &{}= \frac{b^{\hS,\hP}+c^{\hS\Log(\hL)}+\int xyF^{\hS,\Log(\hL)}\big(\d(x,y)\big)}{1+\Delta B^{\Log(\hL),\hP}}; 
\label{eq:b^hS,hPstar}\\
\tilde{c}^{\hS,\hP^\star} &{}= \frac{c^{\hS}+\int xx^{\top }(1+y)F^{\hS,\Log(\hL)}  \big(\d(x,y)\big)}{1+\Delta B^{\Log(L)}}.
\label{eq:c^hShPstar}
\end{align}
\item The set 
$$\hXi _{\ha }=\argmin_{\vt\in\R^d:\vt \hS_{-}=-1}
\{\vt \tilde{c}^{\hS,\hP^\star}\vt^\top-2\vt b^{\hS,\hP^\star}\}$$ 
is non-empty.
\item\label{C:mainhP.1e} For some or, equivalently, any $\hXi_{\ha}$-valued predictable process $\ha$, one has 
\begin{align*}
&\frac{b^{\Log(\hL),\hP}}{1+\Delta B^{\Log(L),\hP}}=-\min_{\vt\in\R^{d}:\vt \hS_{-}=-1}
\{\vt \tilde{c}^{\hS,\hP^\star}\vt^\top-2\vt b^{\hS,\hP^\star}\}
=-\ha\tilde{c}^{\hS,\hP^\star}\ha^{\top }+2\ha b^{\hS,\hP^\star};\\
&-\ha\one_{\zu\tau ,T\zu}\Exp(-(\ha\one_{\zu\tau ,T\zu})\sint \hS)_-\in \adm^{\tau}_{1}(\hS,\hP), 
\end{align*}
for all $[0,T]$-valued stopping times $\tau $.
\end{enumerate}

\item\label{C:mainhP.2} The optimal strategy $\hvp\in\adm_{\hv}(\hS,\hP)$ for the quadratic hedging problems \eqref{eq:3} is given in the feedback form by
\begin{equation}\label{eq:hphi}
\hvp (\hv,\hH)=\hxi +\ha(\hV_{-}-\hvp (\hv,\hH)\hS_{-}),
\end{equation}
where
\begin{equation}\label{eq:hV}
\hV_{t} =\frac{1}{\hL_{t}}\E^{\hP}[\Exp((-\one_{\zu t,T\zu}\ha)\sint \hS)_{T}\hH|\F_{t}],
\end{equation}
$\ha$ is an arbitrary $\hXi_{\ha}$-valued predictable process, and $\hxi$ is an arbitrary predictable process taking values in
$$\hXi_{\smash{\hxi} }=\argmin_{\vt\in\R^d:\vt \hS_{-}
=\hV_{-}}\{\vt\tilde{c}^{\hS,\hP^\star}\vt^{\top }-2\vt\tilde{c}^{\hS\hV,\hP^\star}\}$$ 
with
\begin{equation}\label{eq:201212.1}
\tilde{c}^{\hS\hV,\hP^\star} = \frac{c^{\hS\hV}+\int xz(1+y)F^{\hS,\Log(\hL),\hV,\hP}\big(\d(x,y,z)\big)}{1+\Delta B^{\Log(\hL),\hP}}.
\end{equation}

\item\label{C:mainhP.3} For an arbitrary $\hXi_{\smash{\hxi}}$-valued predictable process $\hxi$, let 
\begin{equation}\label{eq:hateps2hatH}
\hat\eps_t^2(\hH)=\E^{\hP}\left[\one_{\zu t,T\zu}\hL\big(\tilde{c}^{\hV,\hP^\star}
-2\hxi^\top\tilde{c}^{\hS\hV,\hP^\star}+\hxi^\top \tilde{c}^{\hS,\hP^\star}\hxi\big)\sint A_T\,|\,\F_t\right], 
\end{equation}
with
\begin{equation}\label{eq:tchVhPstar}
\tilde{c}^{\hV,\hP^\star}= \frac{c^{\hV}+\int z^2(1+y)F^{\Log(\hL),\hV,\hP}\big(\d(y,z)\big)}{1+\Delta B^{\Log(\hL),\hP}}.
\end{equation} 
Then, the hedging error of the optimal strategy is given by
\begin{equation}\label{eq:mvherrorhP}
\hat\veps^2(\hv,\hH):=\E^{\hP}[(\hvp_T(\hv,\hH)\hS_{T}-\hH)^{2}]=\hL_0(\hv-\hV_0)^2+\hat\eps_0^2(\hH).
\end{equation}
\end{enumerate}
Conversely, but without the assumption that $\cM_2(S,\P)$ contains a strictly positive element, if there exists a semimartingale $\hL=(\hL_t)_{0\leq t\leq T}$ satisfying~\ref{C:mainhP.1}\ref{C:mainhP.1a}--\ref{C:mainhP.1e}, then $\hL$ is the opportunity process and \ref{C:mainhP.2}--\ref{C:mainhP.3} hold.
\end{corollary}

\bpf
Without loss of generality, let $S^1=X$ be the nice numeraire for $S$ under $\P$. Then, $\hS$ is locally square-integrable under $\hP$ by Proposition~\ref{P:ThetaEqual}\ref{ThetaEqual.1bis}.  Recall that $\hS^1=1$ by construction. Theorem~\ref{T:main X=1} applied to the discounted model yields that Theorem~\ref{T:main} holds for the discounted model $(\hS,\hH,\hP)$. Finally, when $\cM_2(S,\P)$ contains a positive deflator $Z$, we have by Proposition~\ref{P:ThetaEqual}\ref{ThetaEqual.4} that $\frac{ZX}{\hZ}>0$ belongs to $\cM_2(\hS,\hP)$, which completes the proof. 
\ep
\section{Optimal hedging without numeraire change.}\label{S:main}
This section contains the main result of the paper, namely a numeraire-invariant explicit representation of optimal
quadratic hedging strategies for the problem \eqref{eq:1new}. Its \emph{proof} relies on the simpler characterization in properly discounted markets, see \cite{cerny.kallsen.07}, which was cast in a symmetric form in Subsection~\ref{SS:hX=1}. From now on, we assume that $S=(X,Y)$ is a locally square-integrable semimartingale.

\subsection{Main results.}\label{SS:5.2}
The first theorem reduces the hedging problem \eqref{eq:1new} to affinely constrained quadratic minimization in $\R^d$
involving the semimartingale characteristics of the price process. Its proof is to be found in Appendix~\ref{S:proof main}. Its conditions and applicability are illustrated in Section~\ref{S:examples}. The extent to which this theorem holds \emph{without} an existence of a nice numeraire is the subject of an ongoing research.
\begin{theorem}[Quadratic hedging without numeraire change]\label{T:main} 
Suppose that $S$ admits a nice numeraire. If $\cM_2(S,\P)$ contains a strictly positive element, the following statements hold.
\begin{enumerate}[(1)]
\item\label{T:main.1} The opportunity process $L$ is the unique bounded semimartingale $L=(L_t)_{0\leq t\leq T}$ such that
\begin{enumerate}[(a)]
\item\label{T:main.1a} $L>0$, $L_->0$, and $L_T=1$.
\item\label{T:main.1b} $\frac{L}{\Exp(B^{\Log(L)})}>0$ is a martingale on $[0,T]$.
\item\label{T:main.1c} $[S,S]$ is $\P^\star$-special for the measure $\P^\star\sim \P$ defined by  
$$\frac{\d\P^\star}{\d\P}=\frac{L_T}{\E[L_0]\Exp(B^{\Log(L)})_T}>0,$$
which implies
\begin{align}
b^{S\star}         &= \frac{b^{S}+c^{S\,\Log(L)}+\int xyF^{S,\Log(L)}\big(\d(x,y)\big)}{1+\Delta B^{\Log(L)}}; \label{eq:b^Sstar}\\
\tilde{c}^{S\star} &= \frac{c^{S}+\int xx^{\top }(1+y)F^{S,\Log(L)}  \big(\d(x,y)\big)}{1+\Delta B^{\Log(L)}}. \label{eq:c^Sstar}
\end{align}
\item\label{T:main.1d} The set 
\begin{equation}\label{eq:min:a}
\Xi _{a}:=\argmin_{\vt\in\R^d:\vt S_{-}=-1}
\{\vt \tilde{c}^{S\star}\vt^\top-2\vt b^{S\star}\}
\end{equation}
is non-empty.
\item\label{T:main.1e} For some or, equivalently, any $\Xi_a$-valued predictable process $a$, one has 
\begin{align}
&\frac{b^{\Log(L)}}{1+\Delta B^{\Log(L)}}=-\min_{\vt\in\R^{d}:\vt S_{-}=-1}
\{\vt \tilde{c}^{S\star}\vt^\top-2\vt b^{S\star}\}\ 
=-a\tilde{c}^{S\star}a^{\top }+2a b^{S\star}, \label{eq:bLog(L)}\\
&-a\one_{\zu\tau ,T\zu}\Exp(-(a\one_{\zu\tau ,T\zu})\sint S)_-\in \adm^{\tau}_{1}(S,\P),\notag
\end{align}
for all $[0,T]$-valued stopping times $\tau $.
\end{enumerate}
\item\label{T:main.2} The optimal strategy $\vp\in\adm_v(S,\P)$ for the quadratic hedging problems \eqref{eq:3} is given in the feedback form by
\begin{equation}\label{eq:mvhst}
\vp (v,H)=\xi +a(V_{-}-\vp (v,H)S_{-}),
\end{equation}
where
\begin{equation}\label{eq:V}
V_{t} =V_{t}(H)=\frac{1}{L_{t}}\E[\Exp((-\one_{\zu t,T\zu}a)\sint S)_{T}H|\F_{t}],
\end{equation}
$a$ is an arbitrary $\Xi_a$-valued predictable  process, and $\xi$ is an arbitrary predictable process taking values in 
\begin{equation}
\Xi_{\xi }:=\argmin_{\vt\in\R^d:\vt S_{-}
=V_{-}}\{\vt\tilde{c}^{S\star}\vt^{\top }-2\vt\tilde{c}^{SV\star}\}\label{eq:min:xi}
\end{equation}
with
\begin{equation}\label{eq:c^SVstar}
\tilde{c}^{SV\star} = \frac{c^{SV}+\int xz(1+y)F^{S,\Log(L),V}\big(\d(x,y,z)\big)}{1+\Delta B^{\Log(L)}}.
\end{equation}

\item\label{T:main.3} 
For an arbitrary $\Xi_\xi$-valued predictable process $\xi$, let
\begin{equation}\label{eq:eps2H}
\eps_t^2(H)=\E[\one_{\zu t,T\zu}L\left(\tilde{c}^{V\star}-2\xi^\top\tilde{c}^{SV\star}
+\xi^\top \tilde{c}^{S\star}\xi\right)\sint A_T\,|\,\F_t], 
\end{equation}
with
\begin{equation}\label{eq:c^Vstar}
\tilde{c}^{V\star}= \frac{c^{V}+\int z^2(1+y)F^{\Log(L),V}\big(\d(y,z)\big)}{1+\Delta B^{\Log(L)}}.
\end{equation}
Then, the hedging error of the optimal strategy is given by
\begin{equation}\label{eq:mvherror}
\veps^2(v,H):=\E\big[\big(\vp_T(v,H)S_{T}-H\big)^{2}\big]=L_0(v-V_0)^2+\eps_0^2(H).
\end{equation}
\end{enumerate}
Conversely, but without the assumption that $\cM_2(S,\P)$ contains a strictly positive element, if there exists a bounded semimartingale $L=(L_t)_{0\leq t\leq T}$ satisfying \ref{T:main.1}\ref{T:main.1a}--\ref{T:main.1e}, then $L$ is the opportunity process and \ref{T:main.2} and \ref{T:main.3} hold.
\end{theorem} 
\begin{remark}[Miscellaneous]
\begin{enumerate}[(1),wide=0pt]
\item Recall that the assumed existence of a positive deflator in $\mathcal M_2(S,\mathsf{P})$ is related to the absence of arbitrage, see Remark~\ref{R:deflators}\ref{R:deflators2}.\smallskip
\item The structure of the solution is reminiscent of the related simpler setup in \cite{cerny.kallsen.07} in that 
the optimal hedging strategy $\vp(v,H)$ consists of two parts. The \emph{pure hedge} $\xi$ 
invests in $S$ in order to reduce the random fluctuations caused by the hypothetical \emph{value process} $V$ of the option $H$.
The second term in \eqref{eq:mvhst}, on the other hand, takes care of the accrued hedging error from the past.
It involves the \emph{adjustment process} $a$ which is -- up to rescaling by current wealth and change of sign -- the optimizer $\vt$
in \eqref{D:op}.\smallskip
\item Identity \eqref{eq:mvhst} yields that $\vp(v,H)=\xi+(V_--W_-)a$, where the wealth process $W=\vp(v,H) S$ solves the stochastic differential equation (SDE) $W=W_-\sint(-a\sint S)+v+(\xi+V_-a)\sint S$. The unique solution to this affine SDE can be represented explicitly, see for example Eberlein and Kallsen \cite[Proposition~3.48]{eberlein.kallsen.19}.\qed
\end{enumerate}\smallskip
\end{remark}

In the second main result of this section, we provide explicit formulae for the sets \eqref{eq:min:a} and \eqref{eq:min:xi} of all minimizers as well as a particular choice of minimizers in terms of semimartingale characteristics and the Moore--Penrose pseudoinverse. The important role of the oblique projector $p\tilde{c}^{S\star}$ appearing in \eqref{eq:a} below is examined in more detail in Appendix~\ref{App:A}.
\begin{theorem}[Explicit expressions]\label{T:explicit}\leavevmode 
\begin{enumerate}[(1)]
\item\label{T:explicit:1} If \eqref{eq:b^Sstar} and \eqref{eq:c^Sstar} are well-defined, then $\Xi_a$ in \ref{T:main.1}\ref{T:main.1d} is non-empty if and only if  $b^{S\star}\in \Ran(\tilde{c}^{S\star})+\Ran(S_-)$, where $\Ran(\cdot)$ denotes the column space of a matrix. In this case, 
$$\Xi _{a} = a + \Null(\tilde{c}^{S\star})\cap \Null(S_{-}^{\top }),$$
where
\begin{equation}\label{eq:a}
a   = (b^{S\star})^\top p  -      \frac{S_{-}^\top}{S_{-}^{\top }S_{-}}(I-\tilde{c}^{S\star}p); \qquad
p = (m\tilde{c}^{S\star}m)^{-1}; \qquad m = I-\frac{S_{-}S_{-}^{\top }}{S_{-}^{\top }S_{-}}.
\end{equation} 
and $\Null(\cdot)$ is the null space of a matrix. Moreover, $a$ is the minimum norm element of $\Xi _{a}$.

\item\label{T:explicit.2} In the setting of Theorem~\ref{T:main}, we have 
$\Xi _{\xi } {}=\xi+\Null(\tilde{c}^{S\star})\cap \Null(S_{-}^{\top })$ with
\begin{equation}\label{eq:xi}
\xi {}= \tilde{c}^{VS\star}p + V_{-}\frac{S_{-}^\top}{S_{-}^{\top }S_{-}}(I-\tilde{c}^{S\star}p)
\end{equation}
and $p$, $m$ as above. Furthermore, $\xi$ is the minimum norm element of $\Xi_{\xi }$.
 
\item\label{T:explicit.3} If $S_{-}\in \Ran(\tilde{c}^{S\star})$, then $a$, $\xi$ in \eqref{eq:a} and \eqref{eq:xi} can alternatively be written as
\begin{align}
a  ={}&(b^{S\star})^\top(\tilde{c}^{S\star})^{-1}-\left(1+(b^{S\star})^\top(\tilde{c}^{S\star})^{-1} S_-\right)
   \frac{S_{-}^{\top }(\tilde{c}^{S\star})^{-1}}{S_{-}^{\top }(\tilde{c}^{S\star})^{-1}S_{-}},  \label{eq: simpl_nice1b} \\ 
\xi ={}&\tilde{c}^{VS\star}(\tilde{c}^{S\star})^{-1}+\left(V_{-}-\tilde{c}^{VS\star}(\tilde{c}^{S\star})^{-1}S_-\right)\frac{S_{-}^{\top }(\tilde{c}^{S\star})^{-1}}{S_{-}^{\top }(\tilde{c}^{S\star})^{-1}S_{-}}.\label{eq: simpl_nice2b}
\end{align}

\item\label{T:explicit.4} If $S_{-}\notin \Ran(\tilde{c}^{S\star})$, let 
$\bar{m}=I-\tilde{c}^{S\star}(\tilde{c}^{S\star})^{-1}$, $\alpha=\frac{S_{-}^{\top }\bar{m}}{S_{-}^{\top }\bar{m}S_{-}}$, and
$r = \alpha b^{S\star}$. Then, $a$, $\xi$ in \eqref{eq:a} and \eqref{eq:xi} can be written as
\begin{align*}
a ={}&(b^{S\star}-rS_{-})^{\top }(\tilde{c}^{S\star})^{-1}-\left(1+(b^{S\star}-rS_{-})^{\top }(\tilde{c}^{S\star})^{-1}S_-\right)\alpha,\\  
\xi ={}&\tilde{c}^{VS\star}(\tilde{c}^{S\star})^{-1}+\left(V_{-}-\tilde{c}^{VS\star}(\tilde{c}^{S\star})^{-1}S_-\right)\alpha.  
\end{align*}
Here $\alpha$ is the instantaneously risk-free fully invested portfolio ($\alpha S_-=1$) and $r$ the instantaneously risk-free rate of return in the sense that the process $r\sint A = \alpha\sint S$ is continuous and satisfies $[\alpha\sint S,\alpha\sint S]=0$.
\end{enumerate}
\end{theorem}
\bpf
By part~\ref{T:constrained q.1} of Theorem \ref{T:constrained q}, $b^{S\star}\in \Ran(\tilde{c}^{S\star})+\Ran(S_-)$ if and only if $\Xi_a$ in \eqref{eq:min:a} is non-empty. Since $\tilde c^{SV\star}\in\Ran(\tilde{c}^{S\star})$, this also yields that $\Xi_{\xi}$ in \eqref{eq:min:xi} non-empty. Formulae \eqref{eq:hatx}, \eqref{eq:hatx2}, and \eqref{eq: modifP} now yield the explicit expressions for $a$ and $\xi$.
\ep
\begin{remark}[Structure condition and the null strategies]\label{R:210421}
The condition 
\begin{equation}\label{eq:SCstar}
b^{S\star}\in \Ran(\tilde{c}^{S\star})+\Ran(S_-)
\end{equation}
in item~\ref{T:explicit:1} is a weak form of an absence of arbitrage under $\P^\star$. In the discounted model, it reduces to the so-called structure condition (SC) introduced in Schweizer \cite{schweizer.95.saa} when viewed under $\hP^\star$, i.e., $b^{\hS,\hP^\star}\in\Ran\big(\tilde{c}^{M^{\hS,\hP^\star},\hP^\star}\big)=\Ran(\tilde c^{\hS,\hP^\star} )$. Condition~\eqref{eq:SCstar} therefore represents an undiscounted form of the classical SC under $\P^\star$. Observe further that the set of $\Null(\tilde{c}^{S\star})\cap \Null(S_{-}^{\top })$-valued strategies coincides with the set of all self-financing strategies with zero wealth, here denoted by $\NullStrategy$. Since the absence of arbitrage as well as the set of null strategies $\NullStrategy$ are measure-invariant, one then obtains an equivalent undiscounted SC in terms of $\P$-characteristics, i.e., $b^{S}\in \Ran(\tilde{c}^{S})+\Ran(S_-)$ and $\NullStrategy$ is equivalently characterized as the set of all $\Null(\tilde{c}^{S})\cap \Null(S_{-}^{\top })$-valued processes. \qed
\end{remark}
\bpf
Under the condition that $b^{S\star}\in \Ran(\tilde{c}^{S\star})+\Ran(S_-)$, we need to argue that an $\R^d$-valued predictable processes $\vp=(\vp_t)_{0\leq t\leq T}$ is a self-financing trading strategy with zero wealth process, that is, $\vp\in\NullStrategy$, if and only if it is valued in $\Null(\tilde{c}^{S\star})\cap\Null(S_-^\top)$.

Suppose first that $\vp\in\NullStrategy$. Then, $\vp S_{-}\equiv 0 $ and $\vp \sint S\equiv 0$. The latter yields that $\langle \vp\sint S\rangle^{\P^\star}=(\vp\, \tilde{c}^{S\star}\, \vp^\top)\sint A\equiv 0$ and therefore $\vp\in\Null(S_-^\top)\cap \Null(\tilde{c}^{S\star})$.

For the converse, assume that $\vp=(\vp_t)_{0\leq t\leq T}$ is valued in $\Null(\tilde{c}^{S\star})\cap\Null(S_-^\top)$. Let $S=S_0+M^{S\star}+B^{S\star}$ be the canonical decomposition of the locally square-integrable and hence special semimartingale $S$ under $\P^\star$. Because $b^{S\star}\in \Ran(\tilde{c}^{S\star})+\Ran(S_-)$, all processes $\vp=(\vp_t)_{0\leq t\leq T}$ valued in $\Null(\tilde{c}^{S\star})\cap\Null(S_-^\top)$ are also valued in $\Null(b^{S\star})$. Then, it follows from the construction of the stochastic integral that $\vp\in L(S)$ with $\vp\sint S\equiv 0$, since $|\vp\sint B^{S\star}|= |\vp\, b^{S\star}|\sint A\equiv0$ and $\langle \vp\sint M^{S\star}\rangle^{\P^\star}\leq \langle \vp\sint S\rangle^{\P^\star}=(\vp\, \tilde{c}^{S\star}\, \vp^\top)\sint A\equiv 0$. Hence, $\vp S_{-}\equiv 0$ and $\vp_0 S_0+\vp\sint S_-\equiv 0$ and $\vp$ is therefore a self-financing strategy with zero wealth process by Proposition \ref{P:sf}.
\ep
\begin{remark}[Universal results and special cases]
The formulae in item~\ref{T:explicit.2} are universal. The special case~\ref{T:explicit.3} is relevant in `discrete-time models.' It always applies when $S$ has no quasi-left-continuous component, i.e., when $S$ can be written as the (not necessarily absolutely convergent) sum of its jumps at predictable stopping times; see \v{C}ern\'y and Ruf \cite[Propositions~3.15 and 4.6]{cerny.ruf.21.bej}. There may or may not be an  instantaneously risk-free asset in item~\ref{T:explicit.3}; this makes no difference to the formulae \eqref{eq: simpl_nice1b}--\eqref{eq: simpl_nice2b}.

The special case~\ref{T:explicit.4} may arise in quasi-left-continuous models, i.e., models where $S$ does not jump at predictable times. It arises only when there is an instantaneously risk-free asset. This is the only occasion where the risk-free asset appears explicitly in Theorem~\ref{T:explicit}. \qed
\end{remark}
\subsection{Verification.}\label{SS:verification}
For applications, we establish the following verification result. This will allow us to identify the opportunity process in models with independent increments, for example. 
\begin{proposition}[Sufficient conditions for the opportunity process]\label{P:verification}
Suppose that $S$ admits a nice numeraire. 
 Let $L=(L_t)_{0\leq t\leq T}$ be a semimartingale such that
 \begin{enumerate}[(1)]
  \item\label{verification.1} $L_T=1$ and $L$ is bounded from above and below by positive constants, which implies that the right-hand sides of \eqref{eq:b^Sstar} and \eqref{eq:c^Sstar} are well-defined.
  \item\label{verification.2} $\Xi_a$ in \eqref{eq:min:a} is non-empty and \eqref{eq:bLog(L)} holds for some $a\in \Xi_a$.
 \end{enumerate}
 Then $L$ satisfies conditions \ref{T:main.1}\ref{T:main.1a}--\ref{T:main.1e} in Theorem~\ref{T:main} and it is therefore the opportunity process.
\end{proposition}

\bpf
Claim~\ref{verification.1} follows from the local square-integrability of $S$ together with the boundedness of jumps of $\Log(L)$. By hypothesis, the joint semimartingale characteristics of $(S,L)$ satisfy
\begin{align}
&b^{\Log(L)}=-\min_{\vt\in\R^{d}:\vt S_{-}=-1}
\{\vt \bar{c}^{S}\vt^\top-2\vt \bar{b}^{S}\}
=-a\bar{c}^{S}a^{\top }+2a \bar{b}^{S} \label{eq: bellLogell}
\end{align}
for some $a\in\bar{\Xi} _{a}=\argmin_{\vt\in\R^d:\vt S_{-}=-1}
\{\vt \bar{c}^{S}\vt^\top-2\vt \bar{b}^{S}\}$, where $\bar{\Xi} _{a}$ is non-empty and
\begin{align*}
\bar{b}^{S}&= b^{S}+c^{S\,\Log(L)}+\int xyF^{S,\Log(L)}\big(\d(x,y)\big); \\
\bar{c}^{S} &= c^{S}+\int xx^{\top }(1+y)F^{S,\Log(L)}  \big(\d(x,y)\big). 
\end{align*}
Without loss of generality, we may take $S^1=X$ to be the nice numeraire, with $S=(X,Y)$. Furthermore, let $\hZ_t=\E[X_T^2|\F_t]$. We will proceed by verifying that the process $\hat{L}=(\hat{L}_t)_{0\leq t\leq T}$ defined by
$$\hat{L}_t=\frac{L_t X^2_t}{\hZ_t},\quad 0\leq t\leq T,$$
is the opportunity process $\hL=(\hL_t)_{0\leq t\leq T}$ in the discounted model with $\hS=S/X=(1,\hY)$.

To this end, we observe that, since $X$ is a nice numeraire and hence
$$0<\underline{\delta}\leq \frac{\hZ_t}{X^2_t}\leq \overline{\delta}<\infty,\quad 0\leq t\leq T,$$ 
and $L$ is bounded above and below by positive constants by asumption, the process $\hat{L}$ is non-negative, bounded from above and below by positive constants. Then, it follows as in the proof of  Theorem~\ref{T:main}\ref{T:main.1} via Lemma \ref{lem: bar to hat} that $L$ satisfies \eqref{eq: bellLogell} if and only if $\hat{L}$ satisfies
\begin{align*}
&b^{\Log(\hat{L}),\hP}=-\min_{\vt\in\R^{d}:\vt \hS_{-}=-1}
\{\vt \hc ^{S}\vt^\top-2\vt \bar{b}^{S}\}
=-\ha \hc ^{S}\ha ^{\top }+2\ha  \hb ^{S}=\ha \hc ^{S}\ha ^{\top }, 
\end{align*}
for any $\ha \in\hXi _{a}=\argmin_{\vt\in\R^d:\vt \hS_{-}=-1}
\{\vt \hc ^{S}\vt^\top-2\vt \hb ^{S}\}$, where $\hXi _{a}=X_-\bar{\Xi} _{a}$ is non-empty and
\begin{align}
\hb ^{\hS}   &{}=b^{\hS,\hP}+c^{\hS\Log(\hat{L})}+\int ylF^{\hS,\Log(\hat{L}),\hP}\big(\d(y,l)\big), \label{eq: hatb0ell}\\
\hc ^{\hS}   &{}=c^{\hS}+\int yy^{\top }(1+l)F^{\hS,\Log(\hat{L}),\hP}\big(\d(y,l)\big).\label{eq: hatc0ell}
\end{align}
Moreover, since $\hS^1\equiv1$, we have that $b^{\hS^1,\hP}=0$, $c^{\hS^1,\hP}=0$ and $F^{\hS^1,\hP}=0$. Therefore, $\ha ^{2:d}=\big((\hb ^{\hS})^{2:d}\big)^\top\big((\hc ^{\hS})^{2:d}\big)^{-1}$ and $\ha ^1=-1-\ha ^{2:d}\hS^{2:d}$, as explained in the proof of Theorem~\ref{T:main X=1}. This also yields that $\ha \hc ^{S}\ha ^{\top }=\big((\hb ^{\hS})^{2:d}\big)^\top\big((\hc ^{\hS})^{2:d}\big)^{-1}(\hb ^{\hS})^{2:d}$.

Because $b^{\Log(\hat{L}),\hP}=\ha \hc ^{S}\ha ^{\top }$ is non-negative, $B^{\Log(\hat{L}),\hP}$ and  $\Exp(B^{\Log(\hat{L}),\hP})$ are non-decreasing. Therefore, $\frac{\hat{L}}{\Exp(B^{\Log(\hat{L}),\hP})}$ is a $\hP$-local martingale that is bounded from above and hence a true $\hP$-martingale. The latter allows to define a measure $\hP^\star\sim\hP$ by $\frac{\d\hP^\star}{\d\hP}=\frac{\hat{L}_T}{\E^{\hP}[\hat{L}_0]\Exp(B^{\Log(\hat{L})})_T}>0.$ Since $\frac{\d\hP^\star}{\d\hP}$ is bounded, $\hS$ remains a locally square-integrable semimartingale under $\hP^{\star}$ with canonical decomposition $\hS=\hS_0+M^{\hS,\hP^\star}+B^{\hS,\hP^\star}$ and
\begin{align*}
b^{\hS,\hP^\star}  &{}= \frac{b^{\hS,\hP}+c^{\hS\Log(\hat{L})}+\int xyF^{\hS,\Log(\hat{L})}\big(\d(x,y)\big)}{1+\Delta B^{\Log(\hat{L}),\hP}}=\frac{\hb ^{\hS}}{1+\Delta B^{\Log(\hat{L}),\hP}}, \\
\tilde{c}^{\hS,\hP^\star} &{}= \frac{c^{\hS}+\int xx^{\top }(1+y)F^{\hS,\Log(\hat{L})}  \big(\d(x,y)\big)}{1+\Delta B^{\Log(\hat{L})}}=\frac{\hc ^{\hS}}{1+\Delta B^{\Log(\hat{L}),\hP}} 
\end{align*}
by Girsanov's theorem. Therefore,
\begin{align*} 
|\ha \sint B^{\hS,\hP^\star}|&=|\ha b^{\hS,\hP^\star}|\sint A=\ha \hc ^{S}\ha ^{\top }\sint A=B^{\Log(\hat{L}),\hP},\\
\int\ha \d\langle M^{\hS,\hP^\star}\rangle^{\hP\star}\ha ^\top&\leq(\ha \tilde{c}^{\hS,\hP^\star}\ha ^\top)\sint A=(\ha \hc ^{S}\ha ^{\top })\sint A=B^{\Log(\hat{L}),\hP},
\end{align*}
are predictable and of finite variation so that $\ha \in L(\hS)$ and the stochastic integral $\ha \sint\hS$ is well defined.

Let $\tau$ be a $[0,T]$-valued stopping time. By applying Yor's formula twice, we have that $\hat{L}\big(\Exp(-\ha \one_{\rrbracket \tau,T\rrbracket}\sint \hS)\big)^2
=\Exp(W),$ where
$$W=\Log(\hat{L})-2(\ha \one_{\rrbracket \tau,T\rrbracket})\sint \hS+[(a\one_{\rrbracket \tau,T\rrbracket})\sint\hS]-2\big[\Log(\hat{L}),(\ha \one_{\rrbracket \tau,T\rrbracket})\sint \hS \big]+\big[\Log(\hat{L}),[(a\one_{\rrbracket \tau,T\rrbracket})\sint\hS ]\big].$$
Thanks to \eqref{eq: bellLogell}, \eqref{eq: hatb0ell} and \eqref{eq: hatc0ell}, we obtain that $B^{W,\hP}=0$. Therefore, $W$ is a local $\hP$-martingale. Hence, $\Exp(W)=\hat{L}\big(\Exp(-\ha \one_{\rrbracket \tau,T\rrbracket}\sint \hS )\big)^2$ is a non-negative, local $\hP$-martingale and therefore a $\hP$-supermartingale. By the optional sampling theorem, this yields that
$$\sup \left\{\E\left[\hat{L}_\sigma\big(\Exp(-\ha \one_{\rrbracket \tau,T\rrbracket}\sint \hS )_\sigma\big)^2\right]~:~0\leq \sigma \leq T \text{ is a stopping time}\right\}<\infty .$$
Since $\hat{L}$ is bounded below by a positive constant, say $\hat{k}>0$, we have, for all $[0,T]$-valued stopping times $\sigma$, that 
$0\leq \big(\Exp(-\ha \one_{\rrbracket \tau,T\rrbracket}\sint \hS )_\sigma\big)^2\leq\frac{\hat{L}_\sigma}{\hat{k}}\big(\Exp(-\ha \one_{\rrbracket \tau,T\rrbracket}\sint \hS )_\sigma\big)^2$. 
Therefore, $\Exp(-\ha \one_{\rrbracket \tau,T\rrbracket}\sint \hS )$ is an $L^2(\hP)$-semimartingale. 
Hence, $\hvt  =-\ha \one_{\rrbracket \tau,T\rrbracket}\Exp(-\ha \one_{\rrbracket \tau,T\rrbracket}\sint \hS )_-$ is a tame trading strategy on $\llbracket \tau,T\rrbracket$ by Corollary \ref{cor: Exp equality}, since $\ha \hS _-=-1$, so that
$$
\hvt = -\ha \one_{\rrbracket \tau,T\rrbracket}\Exp(-\ha \one_{\rrbracket \tau,T\rrbracket}\sint \hS )_-\in\adm^{\tau}_{1}(\hS ,\hP).
$$
This yields that $\hat{L}=(\hat{L}_t)_{0\leq t\leq T}$ satisfies the properties \ref{C:mainhP.1}\ref{C:mainhP.1a}--\ref{C:mainhP.1e} of Corollary \ref{C:mainhP} and is therefore the opportunity process for $(\hS,\hP)$. It then follows as in the proof of Theorem \ref{T:main} that $L=\hat{L}\frac{\hZ}{X^2}$ has the properties \ref{C:mainhP.1}\ref{C:mainhP.1a}--\ref{C:mainhP.1e} of Theorem \ref{T:main}. It is therefore the opportunity process for the undiscounted model $(S,\P)$ and, by the converse implication of Theorem \ref{T:main}, the conclusions \ref{T:main.2}--\ref{T:main.3} of Theorem \ref{T:main} and \ref{T:explicit:1}--\ref{T:explicit.4} of Theorem \ref{T:explicit} hold.
\ep

\section{Efficient frontier.}\label{S:efficient frontier}
We shall call a strategy $\vt^\star\in\adm_v(S,\mathsf{P})$ \emph{weakly efficient} 
(with initial value $v\in\R$)
if $\vt^\star$ minimizes $\Var(\vt_T S_T)$ over all 
$\vt\in\adm_v(S,\P)$ satisfying 
$\E[\vt_T S_T] = \E[\vt^\star_T S_T]$.
Since $\vt^\star$ is weakly efficient (with initial value $v\neq 0$) 
if and only if $\vt^\star/v$ is weakly efficient (with initial value 1),
it suffices to consider the case $v=1$.

From Hansen and Richard~\cite[Lemma 3.3]{hansen.richard.87} it follows that the set of weakly efficient portfolios in 
$\adm_1(S,\P)$ equals
\[\{\vp(1,0)+\lambda\vp(0,1):\lambda\in\R\}.\]
One gleans easily that 
\begin{align*}
\E\big[\big(\vp_T(1,0)S_T&{}+\lambda \vp_T(0,1)S_T\big)^2\big] = \E\big[\big(\vp_T(1,0)S_T\big)^2\big]
+\lambda^2\E\big[\big(\vp_T(0,1)S_T\big)^2\big];\\
\E\big[\vp_T(1,0)S_T\big]&{}=L_0V_0(1);\qquad \E\big[\big(\vp_T(1,0)S_T\big)^2\big]{}=L_0;\\
\E\big[\vp_T(0,1)S_T\big]&{}=\E\big[\big(\vp_T(0,1)S_T\big)^2\big]=1-\veps^2(0,1)=1-L_0V_0^2(1)-\eps_0^2(1),
\end{align*}
see \eqref{eq:mvherror}, where $V(1)$ and $\eps^2(1)$ are given by \eqref{eq:V} and \eqref{eq:eps2H}, respectively. This can be used to determine the weakly efficient frontier.

For any $R\in\{\vt^\star_T S_T:\vt^\star\in\adm_1(S,\P)\text{ weakly efficient}\}$ 
we have 
\begin{equation}\label{eq:ef1}
\E[R^2]=L_0+\big(1-L_0V_0^2(1)-\eps_0^2(1)\big)^{-1}\big( \E[R]-L_0V_0(1)\big) ^{2}.
\end{equation}
Straightforward algebra then gives the equivalent formula \eqref{eq:efintro} in $(\E[R],\Var(R))$-space.
We conclude that the efficient frontier is characterized by three quantities, namely the initial value of the 
opportunity process, the initial value of the optimal hedge of the constant payoff 1 (called \emph{tracking process} in the introduction), and its corresponding hedging error.

\begin{remark}[Simplifications in the presence of a risk-free asset]
When there is a risk-free asset with constant value one, then $V_0(1)=1$ and $\eps_0^2(1)=0$, hence there is only one key quantity, say $L_0$, driving the whole frontier. In a general setting, the value of $L_0$ does not determine $V_0(1)$ or  $\eps_0^2(1)$, hence the task of computing the frontier is roughly three times more demanding compared to the standard case where a risk-free asset of a constant value exists.\qed
\end{remark}
Recall that a strategy $\vt^\star\in\adm_v(S,\P)$ is called \emph{efficient} (with initial value $v\geq0$) if one of the two equivalent statments holds:
\begin{enumerate}
\item $\vt^\star$ maximizes $\E[\vt_T S_T]$ over all 
$\vt\in\adm_v(S,\P)$ satisfying 
$\Var(\vt_T S_T)\leq\Var(\vt^\star_T S_T)$,
\item $\vt^\star$ minimizes $\Var(\vt_T S_T)$ over all 
$\vt \in\adm_v(S,\P)$ satisfying 
$\E[\vt_T S_T]\geq\E[\vt^\star_T S_T]$.
\end{enumerate}
Once again, it suffices to consider the case $v=1$. On a moment's reflection, efficient strategies in $\adm_1(S,\P)$ coincide with the weakly efficient strategies in $\adm_1(S,\P)$ whose mean is above  $\frac{L_0V_0(1)}{\veps^2(0,1)}$. For $V_0(1)>0$ this yields
$\left\{\vt^\star \in \adm_1(S,\P) \text{ efficient }\right\} = \left\{\vp(1,0)+\lambda\vp(0,1)~:~\lambda \geq \frac{L_0V_0(1)}{\veps^2(0,1)}\right\}.$
 
\section{Examples.}\label{S:examples}
We provide two numerical examples, one in discrete and one in continuous time, that illustrate the numeraire-invariant approach to portfolio optimization. We then conclude with a third, theoretical example in an It\^o semimartingale setting. In all three examples we shall assume that logarithmic returns have independent increments. We will now collect some considerations that apply generally in semimartingale models with independent returns.

To begin with, it is convenient to parametrize the trading strategies not in the number of shares~$\vt$, but rather in the dollar amounts invested  in the individual assets 
$$\pi^i=\vt^iS^i_-,\qquad i=1,\ldots,d.$$ 
This gives $\pi=\vt \diag(S_-)$, where $\diag(S_-)$ is the diagonal matrix with $S_-$ on the diagonal and $0$ elsewhere. The formulae of Theorems \ref{T:main} and \ref{T:explicit} remain applicable if one replaces the characteristics of $S$ by those of $\Log(S)=\diag(S_-)^{-1}\sint S$, $\vt$ by $\pi$, and the conditions $\vt S_-=-1$ and $\vt S_-=V_-$ by $\pi \ones=-1$ and $\pi \ones =V_-$, respectively, where $\ones=(1,\ldots,1)^\top\in\mathbb{R}^d$ is a column vector of ones. 

Next, suppose $S>0$, $S_->0$, and that $\Log(S)$ is a  semimartingale with independent increments which is locally-square integrable or equivalently, satisfies $\Log(S)_T\in L^2(\P)$, which in turn is equivalent to $\tilde{c}^{\Log(S)}$ and $\tilde{c}^{\Log(S)}\sint A_T$ being finite. By Proposition~\ref{P:210425} each of the assets $S^i$, $i\in\{1,\ldots,d\}$ is a nice numeraire. Let 
$${\Xi} _{\lambda}=\argmin_{\vt\in\R^d:\vt \ones=-1}\{\vt \tilde{c}^{\Log(S)}\vt^\top-2\vt b^{\Log(S)}\}.$$ 
By part~\ref{T:constrained q.1} of Theorem \ref{T:constrained q}, ${\Xi} _{\lambda}$ is non-empty if and only if the local no-arbitrage condition
\begin{equation}\label{ex:cond:2}
b^{\Log(S)}\in\Ran(\tilde{c}^{\Log(S)})+\Ran(\ones)
\end{equation} holds. Observe that for any $\lambda\in{\Xi} _{\lambda}$ the process $\Exp(-\lambda\sint \Log(S))$ is the wealth of a self-financing strategy and $-\lambda\sint \Log(S)$ is a locally square-integrable semimartingale with independent increments. Then, exactly as in the proof of Proposition~\ref{P:210425}, the process $K=\log(\Exp(B^{-2\lambda\mkern1mu\sint\mkern1mu\Log(S)+[\lambda\mkern1mu\sint\mkern1mu\Log(S)]}))$ is bounded on $[0,T]$.

Proposition~\ref{P:verification} now yields that the deterministic process
$$
L_t
=\exp(K_T-K_t)
=\Exp(B^{-2\lambda\sint\Log(S)+[\lambda\mkern1mu\sint\mkern1mu\Log(S)]})_t/\Exp(B^{-2\lambda\mkern1mu\sint\mkern1mu\Log(S)
+[\lambda\mkern1mu\sint\mkern1mu\Log(S)]})_T$$
is the opportunity process, $\P^\star=\P$, and $a=\lambda$.
This allows us to compute the ingredients of the formulae \eqref{eq:efintro} and \eqref{eq:ef1} of the weakly efficient  mean--variance frontier explicitly. To simplify the notation, we shall write in the rest of this section
$$b^{\Log(S)}\equiv b;\qquad\tilde{c}^{\Log(S)}\equiv c.$$

\begin{example}[Discrete time, IID returns]\label{E:LiNg}
This example illustrates our general machinery on Example 1 in Li and Ng~\cite[Section~7]{li.ng.00}. There are 3 risky assets whose conditional mean rate of return, $\mu$, and variance-covariance matrix, $\Sigma$, read
\begin{equation*}
\mu =\left( 
\begin{array}{c}
0.162 \\ 
0.246 \\ 
0.228
\end{array}
\right) ,\quad
\Sigma =\left( 
\begin{array}{ccc}
146 & 187 & 145 \\ 
187 & 854 & 104 \\ 
145 & 104 & 289
\end{array}
\right)\times 10^{-4}.
\end{equation*}
The model has $T=4$ times steps. For convenience we shall use an activity process with jumps of size $\Delta A_t=1$ at the discrete dates $t\in\{1,\ldots,T\}$. In process notation
$$
b=\mu ,\qquad c=\Sigma +\mu \mu ^{\top}
=\left( 
\begin{array}{ccc}
4.0844 & 5.8552 & 5.1436 \\ 
5.8552 & 14.5916 & 6.6488 \\ 
5.1436 & 6.6488 & 8.0884
\end{array}
\right) \times 10^{-2}.
$$
One then obtains the following values for the key quantities, 
\begin{equation*} m = I-\frac{\ones\ones^\top}{\ones^\top\ones};\qquad p =(mcm)^{\dag}=\left( 
\begin{array}{rrr}
\frac{58\,640\,000}{816\,487} & -\frac{13\,445\,000}{816\,487} & -\frac{45\,195\,000}{816\,487} \\[1.2ex] 
-\frac{13\,445\,000}{816\,487} & \frac{11\,785\,000}{816\,487} & \frac{1\,660\,000}{816\,487} \\[1.2ex] 
-\frac{45\,195\,000}{816\,487} & \frac{1660\,000}{816\,487} & \frac{43\,535\,000}{816\,487}
\end{array}
\right); 
\end{equation*}
\begin{align*}
a ={}& b^{\top }p-\frac{\ones^{\top }}{\ones^{\top }\ones}(I-cp)
\approx \left[ \begin{array}{ccc} -6.9144 & 1.6238 & 4.2907\end{array}\right]; \\[0.5ex]
b^{\top }pb={}&\frac{582\,399}{1632\,974}\approx 0.35665;\qquad 1-ab =\frac{3030\,887}{4082\,435}\approx 0.74242;
\end{align*}
\begin{align*}
&{}1-2ab+aca^{\top}=\frac{14\,\allowbreak 224\,270\,253}{16\,329\,740\,000}\approx 0.87107;\\
L_{t} ={}&(1-2ab+aca^{\top})^{T-t},\qquad
V_{t}(1) =\left( \frac{1-ab}{1-2ab+aca^{\top}}\right) ^{T-t},\qquad t\in \{0,1,\ldots,T\};
\end{align*}
\begin{align*}
L_{0}\approx{}& 0.57571; \qquad L_{0}V_{0}(1)\approx 0.30381;\\
\eps_{0}^{2}(1) ={}&\left( 1-bpb^\top-\frac{(1-ab)^2}{1-2ab+aca^{\top}}\right) \sum_{t=1}^{T}L_t V_t^2(1)\approx 0.024179.
\end{align*}
From \eqref{eq:ef1} the equation for the weakly efficient frontier in the $(\E[R],\E[R^2])$-space reads
\begin{equation*}
\E[R^2]\approx 0.57571+1.2262(\E[R]-0.30381)^{2},
\end{equation*}
while an equivalent formula \eqref{eq:efintro} in $(\E[R],\Var(R))$-space gives
\begin{equation}\label{eq: LiNgError}
\Var(R)\approx 7.5446\times 10^{-2}+0.22625(\E[R]-1.6466)^{2}.
\end{equation}

All numerical values shown here are accurate to the last digit. However, the corresponding result in \cite[p.~403]{li.ng.00} has a small rounding error  in the last digit of the last figure in \eqref{eq: LiNgError}.\qed
\end{example}
\begin{example}[It\^o semimartingale with independent increments]\label{PII}
Suppose $\Log(S)$ is a locally square-integrable It\^o semimartingale with independent increments. Letting $A_t=t$, $t\in[0,T]$, one has by Theorem~\ref{T:explicit}\ref{T:explicit.2}  
\begin{align*}
a   &{}= b^\top p - \frac{\ones^\top}{\ones^{\top }\ones}(I-cp); \qquad 
p = \left(\left(I-\frac{\ones\ones^{\top }}{\ones^{\top }\ones}\right)c
\left(I-\frac{\ones\ones^{\top }}{\ones^{\top }\ones}\right)\right)^{-1}.
\end{align*} 
In the It\^o semimartingale setting 
$$\zeta = \frac{\ones^\top}{\ones^{\top }\ones}(I-cp)$$ is the myopic minimum variance portfolio, i.e., $\zeta$ minimizes the instantaneous variance rate $\pi c\pi^\top -(\pi b)^2\Delta A$ over fully invested portfolios $\pi$ (i.e., $\pi\ones=1$).

By Theorem~\ref{T:constrained q}\ref{T:constrained q.3} we have
$$\zeta c \zeta^\top = aca^\top - (\zeta+a) c (\zeta+a)^\top = aca^\top-b^\top p b\geq 0. $$
Hence, there is a locally risk-free asset with the risk-free rate 
$$r = \zeta b = b^\top p b -ab$$ 
if and only if $aca^\top-b^\top p b = 0$.

Theorem~\ref{T:main} and \cite[Theorem~4.1]{cerny.ruf.23.spa} further yield for all $t\in[0,T]$ that 
\begin{align*}
L_t&{}=\E\left[\Exp\left(-(a\one_{(t,T]})\sint \Log(S)\right)^2_T\right]=\e^{\int_t^T\left(-2a_sb_s+a_s  c_s a_s^\top\right)\d s};\\
L_tV_t(1)&{}=\E\left[\Exp\left(-(a\one_{(t,T]})\sint \Log(S)\right)_T\right]=\e^{-\int_t^Ta_sb_s\d s}.
\end{align*}
Moreover, $\tilde{c}^{V(1)}=0$ and $\tilde{c}^{V(1)\Log(S)}=0$, since $V(1)$ is continuous and deterministic. Item~\ref{T:explicit.2} of Theorem \ref{T:explicit} hence also gives 
\begin{align*}
\xi(1)=  V(1)\zeta;\qquad \xi(1) c\xi(1)^\top{}=V^2(1)\zeta c\zeta^\top;\qquad
\eps^2_0(1){}=\int_0^T L_sV^2_s(1)\zeta_s c_s \zeta_s^\top \d s.
\end{align*}  

Let us illustrate these results with the numerical inputs from \cite[Section~7]{yao.al.14}:
$$ T = 5;\qquad b^\top = [0.2042\ \ 0.5047\ \ 0.1059\ \ 0.0359];\qquad c = \sigma^2,$$
where the half-vectorized form of the symmetric $4\times 4$ matrix $\sigma $ reads
$$\mathrm{vech}(\sigma) =  [\mbox{\small
 1.8385 \ \ 0.3389 \ \ {-0.5712} \ \  0 \ \ 5.8728 \ \  0.8157 \ \ 0.1766 \ \ 1.0503 \ \ {-0.1164} \ \ 0.4604 }].$$
With these inputs one obtains 
$$ \mathrm{vech}(p) \approx [ \mbox{\small 0.3716\ \ {-0.0548} \ \ 0.2913 \ \ {-0.6082} \ \ 0.0490 \ \
 {-0.1221} \ \   0.1279 \ \  0.7933 \ \ {-0.9625} \ \  1.4428} ];$$
$$ a \approx [-0.1172\ \ 0.0852\ \ {-0.3132}\ \ {-0.6548}];\qquad \zeta \approx [0.1745\ \  {-0.0799}\ \ 0.3605\ \ 0.5450];$$ 
\begin{alignat*}{3}
aca^\top&{}\approx 0.08405358;\qquad & \zeta c\zeta^\top &{}\approx 0.06865944;\qquad & ab &{}\approx -0.03762131;\\
L_0&{}\approx 2.21772301;\qquad & L_0V_0(1) &{}\approx 1.20696211;\qquad & \eps_0^2(1)&{} \approx 0.28028620.
\end{alignat*}
One may now chart the weakly efficient frontier using the formulae \eqref{eq:ef1} and \eqref{eq:efintro}, obtaining
\begin{align*}
\E[R^2]\approx{}& 2.21772+15.9127(\E[R]-1.20696)^{2};\\ 
\Var(R)\approx{}& 0.66328+14.9127(\E[R]-1.28790)^{2}.\tag*{\qed}
\end{align*}

\end{example}
\begin{example}[It\^o semimartingales with independent increments: special cases]\ \newline
Suppose $\Log(S)$ is a locally square-integrable It\^o semimartingale with independent increments.
\begin{enumerate}[(A)]
\item In the case $\ones\in\Ran(c)$, Theorem \ref{T:explicit}\ref{T:explicit.3} yields that
\begin{align*}
\zeta = &{} \frac{\ones^{\top }c^{-1}}{\ones^{\top }c^{-1}\ones};\qquad a{}=b^\top c^{-1}-\left(1+b^\top c^{-1} \ones\right)\zeta;\qquad
   \xi(1)=V(1)\zeta;
  \end{align*}
   \begin{alignat*}{2}
   L_t&{}=\e^{\int_t^T\left(\frac{\left(1+b_s^\top c_s^{-1}\ones\right)^2}{\ones^\top c_s^{-1}\ones}-b_s^\top c_s^{-1}b_s\right)\d s};\qquad&
V_t(1)&{}=\e^{-\int_t^T\frac{1+b_s^\top c_s^{-1}\ones}{\ones^\top c_s^{-1}\ones}\d s};\\
L_tV_t(1)&{}=\e^{\int_t^T\left(\frac{\left(1+b_s^\top c_s^{-1}\ones\right)b_s^\top c_s^{-1}\ones}{\ones^\top c_s^{-1}\ones}-b_s^\top c_s^{-1}b_s\right)\d s};\qquad&
L_tV^2_t(1)&{}=\e^{\int_t^T\left(\frac{\left(b_s^\top c_s^{-1}\ones\right)^2-1}{\ones^\top c_s^{-1}\ones}
-b_s^\top c_s^{-1}b_s\right)\d s};
\end{alignat*}
$$
\eps^2_0(1){}=\big(L \xi(1) c \xi(1)^\top\big)\sint A_T=\int_0^T\frac{L_sV^2_s(1)}{\ones^\top c_s^{-1}\ones}\d s.$$
Since we have excluded predictable jump times, there cannot be a risk-free asset, i.e., in the It\^o semimartingale setting the assumption $\ones\in\Ran(c)$ yields
$$ \zeta c \zeta^\top = \frac{1}{\ones^\top c^{-1} \ones}>0.$$

\item If $\ones\notin\Ran(c)$, then there is a locally risk-free portfolio 
$$\zeta=\frac{\ones^\top(I-cc^{-1})}{\ones^\top(I-cc^{-1})\ones}$$ with the risk-free rate of return 
$r=\zeta b$. From \ref{T:explicit.4} of Theorem \ref{T:explicit} we obtain that
\begin{alignat*}{2}
a&{}=(b-r\ones)^\top c^{-1}-\left(1+(b-r\ones)^\top c^{-1} \ones\right)\zeta;\qquad& \xi(1)&{}=V(1)\zeta;\\
L_t&{}=\e^{\int_t^T\left(2r_s-(b_s-r_s\ones)^\top c_s^{-1}(b_s-r_s\ones)\right)\d s};\qquad& V_t(1)&{}= \e^{\int_t^T -r_s \d s};\\
L_tV^2_t(1)&{}=\e^{\int_t^T-(b_s-r_s\ones)^\top c_s^{-1}(b_s-r_s\ones)\d s};\qquad& \eps^2_0(1)&{}=0.\tag*{\qed}
\end{alignat*}
\end{enumerate}  
\end{example}
\newpage
\begin{APPENDICES}
\bookmarksetupnext{level=-1}
\addcontentsline{toc}{section}{Appendices}
\setcounter{appx}{1}
\section{Oblique projections in constrained quadratic optimization.}\label{App:A}
In this appendix, we shall analyse certain oblique projectors arising in finitely-dimensional affinely constrained quadratic optimization that forms part of the quadratic hedging problem in \eqref{eq:min:a} and \eqref{eq:min:xi}. The oblique projectors of interest to us will be of the form 
\begin{equation}\label{eq:UVUV}
E = U(VU)^{-1}V,
\end{equation} 
featuring the Moore--Penrose pseudoinverse here denoted $(\,\cdot\,)^{-1}$. 

We shall proceed in stages. After discussing basic properties of the pseudoinverse, we will observe that in some circumstances the leftmost $U$ or the rightmost $V$ in \eqref{eq:UVUV} can be dropped (Proposition~\ref{P:A+}). We then provide brief overview of  projectors, concluding with the full characterization of the projector $U(VU)^{-1}V$ (Theorem~\ref{T:ObliqueProjector}). As the final preparatory step, Proposition~\ref{P:new2} describes the solution of constrained least squares with oblique projectors in plain view. We then state and prove the main result, Theorem~\ref{T:constrained q}.

Let $m,n\in\N$. For $A\in \Cx^{m\times n}$,
one commonly defines the null space (kernel) and the column space (range) of $A$, respectively, 
$$\Null(A)=\{x\in \Cx^{n}~:~Ax=0\}; \qquad \Ran(A)=\{Ax\in\Cx^m~:~x\in \Cx^{n}\}.$$ 
 If $\X$ is a linear subspace of $\Cx^{m}$, we denote by $\X^{\bot }$ its orthogonal complement, 
$$\X^{\bot }=\{y\in \Cx^{m}~:~y^{\ast}x=0 \text{ for all } x\in \X\}.$$ 
Here $*$ denotes conjugate transpose. For example, $\Ran(A)^\bot=\Null(A^{\ast})$. 

For a given $A\in \Cx^{m\times n}$, consider the following equalities:
\begin{equation}\label{mp}
AXA=A;\qquad XAX ={}X;\qquad AX=X^\ast A^\ast;\qquad XA=A^\ast X^\ast.
\end{equation}
We shall denote by $A^{\dag}$ the Moore--Penrose pseudoinverse of $A$, i.e., the unique matrix $X$ satisfying conditions \eqref{mp} above; see  Ben-Israel and Greville~\cite[Exercise 1.1]{ben-israel.greville.03}. 

The following properties of the Moore--Penrose pseudoinverse will be used routinely.
\begin{proposition}[Properties of the pseudoinverse]\label{P:A+}
The following statements hold:
\begin{enumerate}[(1)]
\item\label{P:A+.1} $(A^{\ast })^{\dag}=(A^{\dag})^{\ast };$
\item\label{P:A+.2} $\Ran(A^{\dag})=\Ran(A^{\ast });$
\item\label{P:A+.3} $U(VU)^{\dag}=(VU)^{\dag}$ if $U$ is an orthogonal projector (i.e., if $U^{2}=U$ and $U^{\ast }=U$);
\item\label{P:A+.4} $(VU)^{\dag}V=(VU)^{\dag}$ if $V$ is an orthogonal projector;
\item\label{P:A+.5} $A^{\dag}=A^{\ast }(AA^{\ast })^{\dag}=(A^{\ast }A)^{\dag}A^{\ast }.$
\end{enumerate}
\end{proposition}

\bpf
\ref{P:A+.1}\quad This follows easily from the conditions \eqref{mp} above.\smallskip

\noindent\ref{P:A+.2}\quad One has 
\begin{align*}
\Ran(A^{\dag}) \supseteq{}&\Ran(A^{\dag}A)\supseteq \Ran(A^{\dag}AA^{\dag})=\Ran(A^{\dag});\\
\Ran(A^{\ast}) \supseteq{}&\Ran(A^{\ast }(A^{\ast })^{\dag})\supseteq \Ran(A^{\ast}(A^{\ast })^{\dag}A^{\ast })=\Ran(A^{\ast }).
\end{align*}
This yields $\Ran(A^{\dag})=\Ran(A^{\dag}A)=\Ran(A^{\ast }(A^{\dag})^{\ast })=\Ran(A^{\ast}(A^{\ast })^{\dag})=\Ran(A^{\ast })$. The second equality follows from the last property in \eqref{mp} for $A^{\dag}$. The third equality follows from part~\ref{P:A+.1} of this proposition. \smallskip

\noindent\ref{P:A+.3}\quad From \ref{P:A+.2}, we have $\Ran(A^{\dag})=\Ran(A^{\ast })$ for all $A$ and therefore 
$$\Ran((VU)^{\dag})=\Ran(U^{\ast }V^{\ast })\subset \Ran(U^{\ast })=\Ran(U).$$ 
Since $U$ is, by assumption, an orthogonal projector with range $\Ran(U)$, the claim follows.\smallskip

\noindent\ref{P:A+.4}\quad By part~\ref{P:A+.3}, $V^{\ast }(U^{\ast }V^{\ast })^{\dag}=(U^{\ast }V^{\ast })^{\dag}$. It now suffices to take the conjugate transpose on both sides.\smallskip

\noindent\ref{P:A+.5}\quad It is straightforward to check that both $X=A^{\ast }(AA^{\ast })^{\dag}$ and $X=(A^{\ast }A)^{\dag}A^{\ast }$ satisfy the conditions \eqref{mp} above. The statement now follows by the uniqueness of $A^{\dag}$.\ep\medskip

We now move to the description of projectors. Two subspaces $L,M\subseteq\Cx^n$ are \emph{complementary} if $L+M=\Cx^n$ and $L\cap M=\{0\}$. For such subspaces, each element $x\in\Cx^n$ has a unique decomposition $x=y+z$ with $y\in L$ and $z\in M$. We then call $y$ the \emph{oblique projection} of $x$ onto $L$ along $M$ and write $y = P_{L,M}x$. Observe that a projector is by necessity idempotent, i.e., $P^2_{L,M}=P_{L,M}$. Conversely, every idempotent matrix $E^2=E$ projects onto $\Ran(E)=L$ along $\Null(E)=M$. The matrix $U(VU)^{-1}V$ is idempotent thanks to the properties of the pseudoinverse.

The orthogonal projector onto $L$ is denoted by $P_{L}=P_{L,L^{\bot }}$. By the properties  of the Moore--Penrose pseudoinverse in \eqref{mp}, $AA^{\dag}$ is idempotent and Hermitian with range $\Ran(A)$, hence $AA^{\dag}=P_{\Ran(A)}$. By Proposition~\ref{P:A+}\ref{P:A+.2} one then has $A^\dag A = P_{\Ran(A^\ast)}$.  The next result describes the projector $U(VU)^{-1}V$.

\begin{theorem}[Properties of the oblique projector $U(VU)^{-1}V$]\label{T:ObliqueProjector}
Let $n,p,q\in\N$. Given two arbitrary matrices $U\in \Cx^{n\times p}$, $V\in \Cx^{q\times n}$, the matrix $E=U(VU)^{\dag}V$ is
idempotent with range and null space given by
\begin{align*}
\Ran(E)& {}=\Ran(UU^{\ast }V^{\ast })=\Ran(UU^{\ast }V^{\ast }V)=\Ran(U)\cap ((UU^{\ast})^{\dag}(\Ran(U)\cap \Null(V)))^{\bot }, \\
\Null(E)&{}=\Null(U^{\ast }V^{\ast }V)=\Null(UU^{\ast }V^{\ast }V)=\Null(V)\oplus (V^{\ast}V)^{\dag}(\Ran(U)+\Null(V))^{\bot }, 
\end{align*}
where `` $\oplus$'' denotes the direct sum.
\end{theorem}
\bpf See \cite[Theorem~3.1]{cerny.09}. \ep

We next proceed with the description of constrained least squares. We shall write $\rank(A)$ for the dimension of $\Ran(A)$; $\rank(A)$ is commonly called the (column) rank of $A$. It is known that the column rank and the row rank coincide, i.e., $r(A)=r(A^\ast)$.

\begin{proposition}[Affinely constrained least squares; {\protect\cite[Corollary~4.3]{cerny.09}}]\label{P:new2}
\!\!\!\!For $k,m,n\in\N$, let $A_{1}\in \Cx^{m\times n}$, $b_{1}\in \Cx^{m}$, $A_{2}\in \Cx^{k\times n}$, and $b_{2}\in \Cx^{k}$ such that $b_2\in\Ran(A_{2})$. Then, the solution of the constrained least squares minimization 
\begin{equation*}
\argmin_{x\in \Cx^{n}}\left\Vert A_{1}x-b_{1}\right\Vert ^{2},\qquad\text{subject to }A_{2}x=b_{2} 
\end{equation*}
is the set 
$ \Xi = \xi +\Null(A_{1})\cap \Null(A_{2})$,
where  
\begin{equation*}
\xi =(A_{1}M)^{\dag}A_{1}A_{1}^{\dag}b_{1}+(I-(A_{1}M)^{\dag}A_{1})A_{2}^{\dag}b_{2}
\end{equation*}
with $M=I-A_{2}^{\dag}A_{2}$. Furthermore, $\xi$ is the element of $\Xi$ with the smallest Euclidean norm.
\end{proposition}

The prominent role of the oblique projector $(A_1M)^\dag A_1$ in the preceding result indicates we are now ready to state and prove the main theorem. In Theorem~\ref{T:constrained q} the key projector is somewhat hidden; it turns out to be $JC$, thanks to the identity $JCJ=J$. Its role in the solution is fully revealed in step~\ref{T:constrained q.4} of the proof, where $JF$ is rewritten as $(JC)C^\dag(CJ)F$. Special cases then follow by suitably reassembling the null space and range of $JC$.

\begin{theorem}[Affinely constrained quadratic optimization]\label{T:constrained q}
For $k,n\in\N$, let $F\in \R^{n}$ and let $C\in \R^{n\times n}$ be a symmetric positive semidefinite matrix. Consider a quadratic form $q(x)=x^{\top}Cx-2x^\top F$ constrained to an affine subspace $\A=\{x\in \R^{n}~:~Ax=b\}$ with $A\in \R^{k\times n}$ and $b\in \R^{k}$ such that $b\in\Ran(A)$.  The following assertions hold:
\begin{enumerate}[(1)]
\item\label{T:constrained q.1}The quadratic form $q$ is bounded from below on $\A$ if and only if 
\begin{equation}\label{eq: admissibleB1}
F\in \Ran(A^{\top })+\Ran(C).  
\end{equation}
\item\label{T:constrained q.2} Provided condition \eqref{eq: admissibleB1} holds, the minimal value of $q(x)$ on $\A$ is attained on the set  
\begin{equation}\label{eq:Xi}
\Xi = \hat{x}+\Null(C)\cap\Null(A),  
\end{equation}
where $J=(MCM)^{\dag}$, $M=(I-A^{\dag}A)$, and
\begin{equation}\label{eq:hatx}
\hat{x}=JF+(I-JC)A^{\dag}b.
\end{equation} 
Furthermore, $\hat{x}$ has the smallest Euclidean norm among all minimizers of $q$ on $\A$. 
\item\label{T:constrained q.3} For any $x\in \Null(A)$ one has $q(\hat{x}-x)=q(\hat{x})+x^{\top }Cx$. Specifically, for $x=JF$ one obtains
\begin{equation*}
q(\hat{x})=q\big((I-JC)A^{\dag}b\big)-F^{\top}JF.
\end{equation*}
\item\label{T:constrained q.4} The following alternative expressions for $\hat{x}$ and $q(\hat{x})$ apply:
\begin{align}
\hat{x} ={}&(I-P_{\bY,\bX})C^{\dag}(I-P_{\bY,\bX}^{\top })F+P_{\bY,\bX}A^{\dag}b; \label{eq:hatx2}\\
q(\hat{x}) ={}&q(P_{\bY,\bX}A^{\dag}b)-F^{\top}(I-P_{\bY,\bX})C^{\dag}(I-P_{\bY,\bX}^{\top })F,\nonumber
\end{align}
with $\bX ={}\Null(A)$ and $\bY ={}(I-CC^{\dag})\Ran(A^{\top })\oplus ^{\bot }C^{\dag}(\Ran(A^{\top})\cap \Ran(C))$, where `` $\oplus ^{\bot }$'' denotes the direct orthogonal sum. Furthermore, 
\begin{equation}\label{eq: modifP}
P_{\bY,\bX}=\left\{ 
\begin{array}{ll}
C^{\dag}A^{\top }(AC^{\dag}A^{\top })^{\dag}A & \qquad\text{for }\Ran(A^{\top })\subset \Ran(C) \medskip\\ 
(I-CC^{\dag})A^{\top }(A(I-CC^{\dag})A^{\top})^{\dag}A & \qquad\text{for } \Ran(A^{\top})\cap \Ran(C) = \{0\}.
\end{array}
\right.
\end{equation}

\end{enumerate}
\end{theorem}

\bpf \ref{T:constrained q.1}\quad We first show the ``only if'' direction. Completion of squares yields
\begin{align*}
q(x) ={}&(x-C^{\dag}F)^{\top }C(x-C^{\dag}F)-F^{\top}C^{\dag}F-2x^{\top}(I-CC^{\dag})F.
\end{align*}
Denote by $y$ the orthogonal projection of $F$ onto $(\Ran(C)+\Ran(A^{\top}))^{\bot }$. Condition \eqref{eq: admissibleB1} is not satisfied if and only if $y\neq 0$. We have $A^{\dag}b+\lambda y\in \A$ for all $\lambda \in \Cx$. Furthermore, 
\begin{align*}
q(A^{\dag}b+\lambda y) ={}&(A^{\dag}b-C^{\dag}F)^{\top }C(A^{\dag}b-C^{\dag}F)-F^{\top}C^{\dag}F \\
&{}-2(A^{\dag}b)^{\top }(I-CC^{\dag})F-2\lambda \|y\|^{2},
\end{align*}
which is unbounded from below as $\lambda \rightarrow +\infty$. Therefore, $q$ is unbounded from below on $\A$ if \eqref{eq: admissibleB1} is not satisfied.

To show the ``if'' direction of \ref{T:constrained q.1}, assume that \eqref{eq: admissibleB1} is met. Consider the decomposition 
\begin{equation*}
F=CJF+(I-CJ)F.
\end{equation*}
By Proposition \ref{P:A+}, 
\begin{equation}\label{eq: CJ}
\begin{split}
JC={}&M(MCM)^{\dag}MC = M((C^{1/2}M)^{\top}C^{1/2}M)^{\dag}MC^{1/2}C^{1/2}\\
={}&M(C^{1/2}M)^{\dag}C^{1/2} = (C^{1/2}M)^{\dag}C^{1/2}. 
\end{split} 
\end{equation}
By Theorem~\ref{T:ObliqueProjector}, $JC=P_{\X,\Y}$ projects onto
\begin{equation}\label{eq: calX2}
\X=\Ran(MC)=\Null(A)\cap (\Null(A)\cap \Null(C))^{\bot }  
\end{equation}
along
\begin{equation}\label{eq: calY2}
\Y=\Null(MC)=\Null(C)\oplus ^{\bot }C^{\dag}(\Null(A)+\Null(C))^{\bot }.
\end{equation}
Thus one obtains that $CJ=(JC)^{\top }=P_{\Y^{\bot },\X^{\bot }}$ with 
\begin{align}
\Y^{\bot } ={}&\Ran(CM)=\Ran(C)\cap (C^{\dag}(\Ran(C)\cap \Ran(A^{\top })))^{\bot };\nonumber\\
\X^{\bot } ={}&\Null(CM)=\Ran(A^{\top })\oplus ^{\bot }(\Ran(A^{\top})+\Ran(C))^{\bot }\label{eq:Xbot}.
\end{align}

Since $\Y^{\bot }\subseteq \Ran(C)$, we have $\Ran(A^{\top })+\Ran(C)\supseteq \Ran(A^{\top })\oplus \Y^{\bot }$, whereby the complementarity of $\Y^{\bot }$ and $\X^{\bot }$ yields 
\begin{equation}\label{eq:Z1}
\Ran(A^{\top })+\Ran(C)=\Ran(A^{\top })\oplus \Y^{\bot }
\end{equation}
as well as 
\begin{equation}\label{eq:Z2} 
\R^n = (\Ran(A^{\top })\oplus \Y^{\bot })\oplus^\bot (\Ran(A^{\top })+\Ran(C))^\bot.
\end{equation}
Since $I-CJ=P_{\X^\bot,\Y^\bot}$, in view of \eqref{eq:Xbot}, \eqref{eq:Z1}, and \eqref{eq:Z2} one obtains 
\begin{equation}\label{eq:(I-CJ)Z}
(I-CJ)(\Ran(A^{\top })+\Ran(C))\subseteq \Ran(A^{\top }).
\end{equation}
Condition \eqref{eq: admissibleB1} now yields
\begin{equation}\label{eq:201220.1}
(I-CJ)F\in \Ran(A^{\top }).  
\end{equation}

Next, observe that $A^\dag A = A^\top (A^\dag)^\top = P_{\Ran(A^\top)}$. By virtue of \eqref{eq:201220.1}, for all $x\in \A$ one has 
$$ x^{\top}(I-CJ)F = x^{\top}A^\top (A^\dag)^\top (I-CJ)F = (A^\dag b)^\top(I-CJ)F.$$
Hence, for all $x\in\A$,
\begin{equation}\label{eq:q}
\begin{split}
q(x) ={}&x^{\top }Cx-2x^{\top }CJF-2x^{\top}(I-CJ)F \\
={}&\|C^{1/2}(x-JF)\|^{2}-F^{\top}JF-2(A^{\dag}b)^{\top }(I-CJ)F
\end{split}
\end{equation}
is bounded from below if \eqref{eq: admissibleB1} holds.\smallskip

\noindent\ref{T:constrained q.2}\quad By \eqref{eq:q}, on $\A$  the minimizers of $q$ coincide with the minimizers of $\|C^{1/2}(x-JF)\|^{2}$. By Proposition~\ref{P:new2} and in view of $\Null(C^{1/2})=\Null(C)$, these minimizers are of the form \eqref{eq:Xi} with 
\begin{align*}
\hat{x} ={}&  (C^{1/2}M)^{\dag}C^{1/2}(C^{1/2})^{\dag}C^{1/2}JF+(I-(C^{1/2}M)^{\dag}C^{1/2})A^{\dag}b
={} JCJF+(I-JC)A^{\dag}b,
\end{align*}
where the second equality follows from \eqref{eq: CJ}. Finally, the identity $J=JCJ$ yields \eqref{eq:hatx}.\smallskip

\noindent\ref{T:constrained q.3}\quad For any $x\in \Null(A)$ 
\begin{align*}
q(\hat{x}-x) ={}&q(\hat{x})+x^{\top }Cx-2x^{\top }(C\hat{x}+F) \\
={}&q(\hat{x})+x^{\top }Cx-2x^{\top }(I-CJ)(CA^{\dag}b+F)  
={}q(\hat{x})+x^{\top }Cx, 
\end{align*}
where the last equality follows from \eqref{eq:(I-CJ)Z}. The rest of the claim follows easily.\smallskip

\noindent\ref{T:constrained q.4}\quad First, note that $\hat{x} =JCJF+(I-JC)A^{\dag}b =JCC^{\dag}CJF+(I-JC)A^{\dag}b$. 
We have shown previously $JC=P_{\X,\Y}$ with $\X,\Y$ given in \eqref{eq: calX2} and \eqref{eq: calY2}. By virtue of Lemma~\ref{L:orth_decomp}, $\Null(C)$ has the following orthogonal decomposition
\begin{equation*}
\Null(C)=(\Null(C)\cap \Null(A))\oplus ^{\bot }(I-CC^{\dag})\Ran(A^{\top }).
\end{equation*}

Consider now the following complementary subspaces of $\R^n$ obtained by ``moving'' the subspace $\Null(C)\cap \Null(A)$ from $\Y$ to $\X$,
\begin{align*}
\bX ={}&\X\oplus^\bot(\Null(C)\cap \Null(A)) = \Null(A);\qquad
\bY =(I-CC^{\dag})\Ran(A^{\top })\oplus ^{\bot }C^{\dag}(\Ran(A^{\top})\cap \Ran(C)).
\end{align*}
By construction, $JCz=P_{\X,\Y}z=P_{\bX,\bY}\,z$ for all $z\in \Ran(A^\top)+\Ran(C)$.
This means one can replace $JC$ with $P_{\bX,\bY}$ in the formula for $\hat{x}$,
\begin{align*}
\hat{x} ={}&JCC^{\dag}CJF+(I-JC)A^{\dag}b  
=(I-P_{\bY,\bX})C^{\dag}(I-P_{\bY,\bX}^{\top })F+P_{\bY,\bX}A^{\dag}b.
\end{align*}
Since $\bX$ and $\bY$ are complementary, Theorem \ref{T:ObliqueProjector} with $U=C^{\dag}A^\top$ resp. $U = (1-C^{\dag}C)A^\top$ and $V=A$  yields \eqref{eq: modifP}. The alternative formula for $q(\hat{x})$ follows similarly from item \ref{T:constrained q.3}.
\ep

\begin{lemma}[Orthogonal complement of $L\cap M$ in $L$]
\label{L:orth_decomp}For any two subspaces $L$, $M$ of\/ $\Cx^n$, $n\in\N$, the following equalities hold:
\begin{equation*}
L= (L\cap M) \oplus ^{\bot }(L\cap (L\cap M)^{\bot })=(L\cap M)\oplus ^{\bot}P_{L}M^{\bot }.
\end{equation*}
\end{lemma}

\bpf
The first equality is obvious. The second equality follows from 
\begin{align*}
P_{L}M^{\bot }={}&\Ran(P_{L}(I-P_{M}))=\Null((I-P_{M})P_{L})^{\bot }
=(L^{\bot }\oplus^\bot (L\cap M))^{\bot }=L\cap (L\cap M)^{\bot }. \tag*{\qed}
\end{align*}
\newpage
\addtocounter{appx}{1}
\section{Proof of Theorem~\ref{T:main X=1}.}\label{S:proof X=1}
We write $x^{2:d}=(x_2,\ldots,x_d)$ for $x=(x_1,\ldots,x_d)\in\R^d$ and $c^{2:d}=(c_{ij})_{2\leq i,j\leq d}$ for $c=(c_{ij})_{1\leq i,j\leq d}\in\R^{d\times d}$. For any $[0,T]$-valued stopping time $\tau$ let 
$$\adm^{\CK,\tau}_{0}({Y},\P)=\left\{\eta\in\adm^{\CK}({Y},\P)~:~\text{$\eta=0$ on $\llbracket0,\tau\zu$} \right\}.$$
Since
$$
{L}_t=\underset{{\vt}\in\adm^{t}_{1}({S},\P)}{\essinf}\E\left[({\vt}_T{S}_T)^2\big|\F_t\right]\\
=\underset{\eta\in\adm^{\CK,t}_{0}({Y},\P)}{\essinf}\E\left[(1-\eta\sint{Y}_T)^2\big|\F_t\right]$$
by Proposition \ref{P:strategies}, we have that the current definition of the opportunity process in \eqref{D:op} coincides with that given in  \cite[Definition~3.3]{cerny.kallsen.07} and \cite[Proposition~6.1]{czichowsky.schweizer.13}.

By \cite[Lemma~3.17~and~Theorem~3.25]{cerny.kallsen.07}, it follows, under the assumption that $Y$ admits an equivalent local martingale measure with square integrable density, that the opportunity process ${L}$ is the unique semimartingale ${L}=({L}_t)_{0\leq t\leq T}$ such that the following statements hold.
\begin{enumerate}[(A),wide = 0pt]
\item\label{CKL.a} ${L}>0$, ${L}_->0$, ${L}_T=1$, and $L$ is bounded. (It then follows from item~\ref{CKL.b} below that $L$ and $L_-$ are in fact $(0,1]$-valued. Indeed, because ${L}>0$ and ${L}_->0$ by \ref{CKL.a}, we have that $1+\Delta B^{\Log(L)}>0$. Hence, $b^{{L}}\geq0$ by \ref{CKL.b} so that $L$ is a non-negative submartingale with ${L}_T=1$ and therefore bounded above by $1$.)
\item\label{CKL.b} The joint $\P$-semimartingale characteristics of $({Y},{L})$ solve the equation
\begin{equation}
\frac{b^{{L}}}{1+\Delta B^{\Log(L)}}=-{L}_-\min_{\vt^{2:d}\in\R^{d-1}}
\{\vt^{2:d}\tilde c^{Y\star}(\vt^{2:d})^\top-2\vt^{2:d} b^{Y\star}\}
={L}_-(b^{Y\star})^\top\big(\tilde c^{Y\star}\big)^{-1}b^{Y\star},\label{eq: hatbl}
\end{equation}
where  
\begin{align*}
b^{Y\star}    ={}&\frac{b^{{Y}}+c^{{Y}\Log({L})}+\int ylF^{{Y}\Log({L})}\big(\d(y,l)\big)}{1+\Delta B^{\Log(L)}},\\
\tilde c^{Y\star} ={}&\frac{c^{{Y}}+\int yy^{\top }(1+l)F^{{Y}\Log({L})}\big(\d(y,l)\big)}{1+\Delta B^{\Log(L)}}.
\end{align*}
In particular, \eqref{eq: hatbl} implies that $b^{Y\star}$ is in the range (column space) of $\tilde{c}^{Y\star}$.
\item\label{CKL.c} For ${a}^{2:d}=(b^{Y\star})^\top(\tilde{c}^{Y\star})^{-1}$, we have that
\begin{equation*}\label{eq: lambda hat admissible}
-({\lambda}^{(\tau)})^{2:d}=-{a}^{2:d}\one_{\zu\tau,T\zu}
\Exp\big((-{a}^{2:d}\one_{\zu\tau ,T\zu})\sint {Y}\big)_-\in \adm^{\CK}({Y},\P)
\end{equation*}
holds for any $[0,T]$-valued stopping time $\tau$. In particular, ${a}^{2:d}$ meets the conditions of the adjustment process from \cite[Definition~3.8]{cerny.kallsen.07}.
\end{enumerate}

From Lemma 3.27 and Definition 3.28 in \cite{cerny.kallsen.07}, it follows that
\begin{align*}
{L}_t&=\underset{\eta\in\adm^{\CK,t}_{0}({Y},\P)}{\essinf}\E\left[(1-\eta\sint{Y}_T)^2\big|\F_t\right]\\
&=\E\left[\big(1-({\lambda}^{(t)})^{2:d}\sint{Y}_T\big)^2\big|\F_t\right]
=\E\left[\Big(\Exp((-{a}^{2:d}\one_{\zu t ,T\zu})\sint {Y})_T\Big)^2\big|\F_t\right]
\end{align*}
and hence, for any $[0,T]$-valued stopping time $\tau$, the process $({\lambda}^{(\tau)})^{2:d}$ solves the minimization problem 
\begin{equation}\label{ap}
\E\left[(1-\eta\sint{Y}_T)^2\big|\F_\tau\right]\to\min_{\eta\in\adm^{\CK,\tau}_{0}({Y},\P)}!
\end{equation}

Conversely, it follows from part (3) of Proposition 6.1 of \cite{czichowsky.schweizer.13} that a process $L=(L_t)_{0\leq t\leq T}$ with the properties \ref{CKL.a}--\ref{CKL.c} is the opportunity process and, for any $[0,T]$-valued stopping time $\tau$, the strategy $({\lambda}^{(\tau)})^{2:d}={a}^{2:d}\one_{\zu\tau,T\zu}\Exp\big((-{a}^{2:d}\one_{\zu\tau ,T\zu})\sint {Y}\big)_-$ solves the minimization problem \eqref{ap}. Note that this conclusion does not require the existence of an equivalent martingale measure for $Y$, but is, by Theorem 6.2 of \cite{czichowsky.schweizer.13}, equivalent to the weaker condition that $Y$ is a $\big(\Exp(N),L\big)$-martingale for a suitable locally square-integrable local $\P$-martingale $N$ such that $\big(\Exp(N),L\big)$ is regular and square-integrable in the sense of Definitions 2.7 and 2.11 of \cite{czichowsky.schweizer.13}, respectively. Then, Theorem 2.16 of \cite{czichowsky.schweizer.13} yields that the space of stochastic integrals $\{v+\eta\sint Y_T~:~\eta\in\adm^{\CK}_{0}({Y},\P) \}$ is closed in $L^2(\P)$ for any $v \in L^2(\F_{0},\P)$ and hence a solution to the quadratic hedging problem exists and the conclusions of \cite{cerny.kallsen.07} hold. 

By the self-financing condition and the definition of the stochastic exponential, we have on $\zu \tau,T\zu$ that
\begin{align*}
-({\lambda}^{(\tau)})^1&{}=\Exp\big((-{a}^{2:d}\one_{\zu \tau ,T\zu})\sint {Y}\big)_{-}
-\Exp\big((-{a}^{2:d}\one_{\zu \tau ,T\zu})\sint {Y}\big)_{-}(-\one_{\zu \tau ,T\zu}{a}^{2:d}) {Y}_{-}\\
&{}=\Exp\big((-{a}^{2:d}\one_{\zu \tau ,T\zu})\sint {Y}\big)_{-}(1+\one_{\zu \tau ,T\zu}{a}^{2:d}{Y}_{-}) .
\end{align*}
Recall that $Y=S^{2:d}$. Hence, setting ${a}^1=-1-{a}^{2:d} {S}^{2:d}_-$ yields that 
$$-{\lambda}^{(\tau)}=-{a}\one_{\zu \tau ,T\zu}\Exp\big((-{a}\one_{\zu \tau ,T\zu})\sint {S}\big)_{-}$$
is a self-financing trading strategy on $\zu \tau,T\zu$ starting with initial wealth $1$ at time $\tau$ such that 
$-{\lambda}^{(\tau)}S=\Exp\big((-{a}\one_{\zu \tau ,T\zu})\sint {S}\big)$. By Proposition \ref{P:strategies}, 
$$-{\lambda}^{(\tau)}=-{a}\one_{\zu \tau ,T\zu}\Exp\big((-{a}\one_{\zu \tau ,T\zu})\sint {S}\big)_{-}\in\adm^{\tau}_{1}({S},\P)$$
if and only if
$$-({\lambda}^{(\tau)})^{2:d}=-{a}^{2:d}\one_{\zu\tau,T\zu}
\Exp\big((-{a}^{2:d}\one_{\zu\tau ,T\zu})\sint {Y}\big)_-\in\adm^{\CK}({Y},\P).$$

Next, we observe that ${S}^1=X=1$ has $b^{X}=0$, $c^{X}=0$ and $F^{X}=0$ and hence 
\begin{equation}\label{eq:0border}
{b}^{X\star}=0;\qquad \tilde{c}^{X\star}=0;\qquad \tilde{c}^{XY\star}=0.
\end{equation}
Therefore, since ${a}{S}_-=-1$ and
$$b^{\Log({L})} = (b^{Y\star})^\top\big(\tilde{c}^{Y\star}\big)^{-1}b^{Y\star}=-\min_{\vt^{2:d}\in\R^{d-1}}\{\vt^{2:d} \tilde{c}^{Y\star}(\vt^{2:d})^\top-2\vt^{2:d} b^{Y\star}\}$$
by the properties of the Moore--Penrose pseudoinverse, we obtain that
$${a}\in\Xi_a =\argmin_{\vt\in\R^{d}:\vt {S}_{-}=-1}\{\vt \tilde{c}^{S\star}\vt^\top-2\vt {b}^{S\star}\}$$
and that 
$$b^{\Log({L})}=-\min_{\vt\in\R^{d}:\vt {S}_{-}=-1}
\{\vt^\top \tilde{c}^{S\star}\vt-2\vt^\top {b}^{S\star}\}
=-{a}^{\top }\tilde{c}^{S\star}{a}+2{a}^{\top }{b}^{S\star}.$$

Moreover, because $S^1=1$, we have that
$\Exp\big((-{a}\one_{\zu \tau ,T\zu})\sint {S}\big)=\Exp\big((-{a}^{2:d}\one_{\zu \tau ,T\zu})\sint {S}^{2:d}\big)$
for all $[0,T]$-valued stopping times $\tau$ and hence
$$L_t=\E\left[\left(\Exp\big((-{a}\one_{\zu t ,T\zu})\sint {S}\big)_T\right)^2\big|\F_t\right].$$

Since ${L}$ and ${L}_-$ are $(0,1]$-valued, we have that $\Log(L)$ is a special semimartingale and $\frac{L}{\Exp(B^{\Log})}$ is a martingale on $[0,T]$. As in \cite[Lemma~3.15]{cerny.kallsen.07}, this allows the definition of the opportunity neutral measure $\P^{\star}\sim\P$ by setting 
$$ \frac{\d\P^\star}{\d\P}=\frac{L_T}{L_0\Exp(B^{\Log(L)})_T}.$$ 
By \ref{CKL.b}, $b^{Y\star}$ is in the range of $\tilde{c}^{Y\star}$ and therefore $b^{S\star}$ is in the range of $\tilde{c}^{S\star}$ by \eqref{eq:0border}. Consequently, any predictable process $\psi$ valued in the kernel of $\tilde{c}^{Y\star}$ with $\psi S_-=0$ is in $L(S)$ with $\psi\sint S=0$. This yields that $\tilde{a}\sint S={a}\sint S$ for any $\tilde{a}\in{\Xi} _{{a}}$.

To sum up, in the special case $S=(1,Y)$ a stochastic process ${L}$ is the opportunity process in the sense of \eqref{D:op} if and only if it satisfies conditions \ref{CKL.a}--\ref{CKL.c} above. Conditions \ref{CKL.a}--\ref{CKL.c} are equivalent to conditions \ref{T:main.1a}, \ref{T:main.1c} and \ref{T:main.1e} of Theorem~\ref{T:main}. Properties \ref{T:main.1b} and \ref{T:main.1d} of Theorem~\ref{T:main} then follow from \ref{CKL.b}. We have therefore shown that a stochastic process ${L}$ is the opportunity process if and only if it satisfies the conditions \ref{T:main.1}\ref{T:main.1a}--\ref{T:main.1e} of Theorem~\ref{T:main}. 

In order to determine the quadratic hedging strategy ${\vp}({v},{H})$ for the contingent claim ${H}$ with initial capital ${v}\in\R$ under $\P$, we need to characterize the mean value process $V=(V_t)_{0\leq t\leq T}$ and the pure hedge coefficient $\xi=(\xi)_{0\leq t\leq T}$.

By \cite[Definition~4.2]{cerny.kallsen.07}, the mean value process $ V $ is given by
$$ V _t=\E\left[\Exp(\one_{\zu t,T \zu }\sint {N})_T H ~|~\F_t\right],\quad0\leq t\leq T,$$
where
$${N}=\Log(L)-{a}^{2:d}\sint Y-[{a}^{2:d}\sint Y,\Log(L)]=\Log(L)-{a}\sint S -[{a}\sint S ,\Log(L)]$$
is the variance-optimal logarithm process. Therefore, it follows from Yor's formula and $L_T=1$ that
 $V _t={}\E\left[\Exp((-\one_{\zu t,T \zu }{a})\sint {S})_T\Exp(\one_{\zu t,T \zu }\sint \Log(L))_T H ~|~\F_t\right]
={}\frac{1}{L_t}\E\left[\Exp((-\one_{\zu t,T \zu }{a})\sint {S})_T H ~|~\F_t\right]$,
which yields \eqref{eq:V}.

Recall $Y=S ^{2:d}$. Combining \cite[Definition~4.6]{cerny.kallsen.07} with 
$$\langle Y,Y\rangle^{\P^\star}=\tilde{c}^{Y\star}\sint A;\qquad \langle Y,V\rangle=\tilde{c}^{YV\star}\sint A,$$ 
where 
$$\tilde{c}^{SV\star}=\frac{c^{SV}+\int xz(1+y)F^{S,\Log(L),V}(\d(x,y,z))}{1+\Delta B^{\Log(L)}}$$ 
by formula (4.16) of \cite{cerny.kallsen.07}, gives 
$$\xi^{2:d}=\tilde{c}^{VS\star}(\tilde{c}^{S\star})^{-1}$$ 
for the pure hedge coefficient $\xi$. By Theorem~4.10 of \cite{cerny.kallsen.07}, we then have
$$\vp^{2:d}_t({v}, H )=\xi^{2:d}_t+\left( V _{t-}-{v}-\vp^{2:d}({v}, H )\sint S ^{2:d}_{t-}\right){a}^{2:d}_t,\quad 0\leq t\leq T,$$
for the optimal hedging strategy $\vp(v, H )$.

From the self-financing condition, we obtain
\begin{align*}
\vp_t^{1}({v}, H )&={v}+\vp^{2:d}({v}, H )\sint S ^{2:d}_{t-}-\vp^{2:d}_t({v}, H )S ^{2:d}_{t-}\\
&={v}+\vp^{2:d}({v}, H )\sint S ^{2:d}_{t-}
-\left(\xi^{2:d}_t+\left( V _{t-}-{v}-\vp^{2:d}({v}, H )\sint S ^{2:d}_{t-}\right){a}^{2:d}_t\right)S^{2:d}_{t-}\\
&=( V _{t-}-\xi^{2:d}_tS ^{2:d}_{t-})
+\left( V _{t-}-{v}-\vp^{2:d}(c, H )\sint S ^{2:d}_{t-}\right)(-1-{a}^{2:d}_tS ^{2:d}_{t-}),\quad 0\leq t\leq T.
\end{align*}
Since ${a}^1=-1-{a}^{2:d}S ^{2:d}_{-}$, setting $\xi^1= V _--\xi^{2:d}S ^{2:d}_{-}$ gives
\begin{equation*}
\vp ({v},H)=\xi +( V _{-}-\vp ({v}, H )S_{-}){a}.
\end{equation*}

Moreover, by the properties of the Moore--Penrose pseudoinverse, $\xi^{2:d}=\tilde{c}^{VY\star}(\tilde{c}^{Y\star})^{-1}$ is the solution to 
$$\vt^{2:d}\tilde{c}^{Y\star}(\vt^{2:d})^\T-2\vt^{2:d}\tilde{c}^{YV\star}\to\min_{\vt^{2:d}\in\R^{d-1}}!$$
Since $\tilde{c}^{VX\star}=\bar{c}^{ V S ^1}=0$, $\tilde{c}^{XY\star}=(\tilde{c}^{YX\star})^\top=0$, and $\xi S _-= V _-$, we have that
$\xi$ solves the constrained minimisation problem
$\vt\tilde{c}^{S\star}\vt^\top-2\vt\tilde{c}^{SV\star}\to\min_{\vt\in\R^{d}:\vt S _-= V _-}!$
and
\begin{align*}
\min_{\vt\in\R^{d}:\vt S _-= V _-}\{\vt\tilde{c}^{S\star}\vt^\top-2\vt\tilde{c}^{SV\star}\}&{}=
\min_{\vt^{2:d}\in\R^{d-1}}\{\vt^{2:d}\tilde{c}^{Y\star}(\vt^{2:d})^\top-2\vt^{2:d}\tilde{c}^{YV\star}\}
{}=-\tilde{c}^{VY\star}(\tilde{c}^{Y\star})^{-1}\tilde{c}^{YV\star}.
\end{align*}
For any self-financing trading strategy, $v+\vp(v, H )\sint  S _{-}=\vp(v, H )S _{-}$ by Proposition~\ref{P:sf}\ref{P:sf.ii} and hence \eqref{eq:mvhst} holds.

By Theorem 4.12 of \cite{cerny.kallsen.07}, the mean squared hedging error of the optimal strategy is given by
$L_0({v}- V _0)^2+\E[L\sint\langle  V -\xi^{2:d}\sint Y\rangle^{\P^{\star}}_T]$.
Because $ Y= S ^{2:d}$, $\xi^{2:d}=\tilde{c}^{VY\star}(\tilde{c}^{Y\star})^{-1}$, and $\langle V,V\rangle^{\P^\star}=\tilde{c}^{V\star }\sint A$, where
$$\tilde{c}^{V\star}=\frac{c^{ V }+\int z^2(1+y)F^{\Log(L), V}\big(\d(y,z)\big)}{1+\Delta B^{\Log(L)}},$$
by formula~(4.17) of \cite{cerny.kallsen.07}, we have that 
\begin{align*}
\langle  V -\xi^{2:d}\sint Y\rangle^{\P^{\star}}&{}=\langle  V,V \rangle^{\P^{\star}}-\xi^{2:d}\sint\langle Y,V\rangle^{\P^\star}
=\left(\tilde{c}^{V\star}-\tilde{c}^{VY\star}(\tilde{c}^{Y\star})^{-1}\tilde{c}^{YV\star}\right)\sint A\\
={}&\Big(\tilde{c}^{V\star}+\min_{\vt\in\R^{d}:\vt S _-= V _-}\{\vt\tilde{c}^{S\star}\vt^\top-2\vt\tilde{c}^{SV\star}\}\Big)\sint A
{}=\big(\tilde{c}^{V\star}-2\xi\tilde{c}^{SV\star}+\xi\tilde{c}^{S\star}\xi^\top\big)\sint A.
\end{align*}
Therefore
$\E\big[L\sint\langle  V -\xi^{2:d}\sint Y\rangle^{\P^{\star}}\big]
=\E\left[L\big(\tilde{c}^{V\star}-2\xi\tilde{c}^{SV\star}+\xi\tilde{c}^{S\star}\xi^\top\big)\sint A_T\right]$, 
which proves \eqref{eq:mvherror}.

This completes the proof of the properties \ref{T:main.1}--\ref{T:main.3} in Theorem~\ref{T:main} for $S=(1,Y)$.
\setlength{\abovedisplayskip}{7pt}
\setlength{\belowdisplayskip}{7pt}
\newpage
\addtocounter{appx}{1}
\section{Proof of Theorem~\ref{T:main}.}\label{S:proof main}
Without loss of generality, let $S^1=X$ be a nice numeraire. 
Observe that \eqref{eq:4} holds thanks to Proposition~\ref{P:ThetaEqual}\ref{ThetaEqual.1}. Corollary~\ref{C:mainhP} yields the solution for the discounted model $(\hS,\hH,\hP)$. We shall now exploit the equalities \eqref{eq:3} and \eqref{eq:3.5} to obtain the solution for the undiscounted model $(S,H,\P)$. 

\ref{T:main.1} Because $X$ is a nice numeraire, we have that $\tame^\tau(\hS,\hP)=\tame^\tau(S,\P)$ by Proposition \ref{P:ThetaEqual}. Since $X_\tau\in(0,\infty)$ is $\F_\tau$-measurable, this implies $\tilde{\vt}\in\tame^\tau_{1}(\hS,\hP)$ if and only if  $\tilde{\vt}/X_\tau\in\tame^\tau_{1}(S,\P)$
from Proposition~\ref{P:sf}\ref{P:sf.iii} and hence
\begin{equation}\label{eq:210421}
\frac{1}{X_\tau}\adm^{\tau}_{1}(\hS,\hP)=\adm^{\tau}_{1}(S,\P).
\end{equation}

Combining this with the definition of $\hP$ via $\frac{\d\hP}{\d\P}=\frac{X_T^2}{\E[X_T^2]}$ and $\hZ_t=\E[X_T^2|\F_t]$, we obtain using Bayes' formula and $X_T\tilde{\vt}_T\hS_T=\tilde{\vt}_TS_T$ for $\tilde{\vt}\in\adm^{t}_{1}(\hS,\hP)$ that 
\begin{align*}
\hL_t&{}=\underset{\tilde{\vt}\in\adm^{t}_{1}(\hS,\hP)}{\essinf}\E^{\hP}\left[(\tilde{\vt}_T\hS_T)^2\big|\F_t\right]\\
&{}=\underset{\tilde{\vt}\in\adm^{t}_{1}(\hS,\hP)}{\essinf}
\left(\frac{X^2_t}{\hZ_t}\E\left[\left.\left(\frac{\tilde{\vt}_T}{X_t}S_T\right)^2\,\right|\F_t\right]\right)
{}=\frac{X^2_t}{\hZ_t}\underset{\tilde{\psi}\in\adm^{t}_{1}(S,\P)}{\essinf}
\E\left[(\tilde{\psi}_TS_T)^2|\F_t\right]=\frac{X^2_t}{\hZ_t}L_t.
\end{align*}

Therefore, under the assumption that $X$ is a nice numeraire, $L$ is the opportunity process for the price process $S$ under $\P$ if and only if $\hL=\frac{X^2}{\hZ}L$ is the opportunity process for the price process $\hS$ under $\hP$. Moreover, by the invariance of semimartingales under equivalent changes of measure and because of the uniform bounds 
$0<\underline{\delta}\leq\frac{X^2}{\hZ}\leq\overline{\delta}<\infty$ in \eqref{eq:vnice}, we have that $L$ is a bounded semimartingale $L=(L_t)_{0\leq t\leq T}$ with $L,L_->0$ and $L_T=1$ if and only if $\hL=\frac{X^2}{\hZ}L$ is a bounded $\hP$-semimartingale $\hL=(\hL_t)_{0\leq t\leq T}$ with $\hL,\hL_->0$ and $\hL_T=1$.

In addition, under the assumption that $X$ is a nice numeraire, we have that 
$$(1,0\ldots,0)\one_{\zu \tau,T\zu}\in\adm^{\tau}_{1}(S,\P).$$ 
Hence, $L\leq\frac{\hZ}{X^2}$ and $\hL\leq 1$. For $L=\frac{\hZ}{X^2}\hL$ bounded, the expressions in \eqref{eq:b^Sstar} and \eqref{eq:c^Sstar} are well defined because $S$ is locally square-integrable under $\P$.

It follows from \eqref{eq:b^hS,hPstar}, \eqref{eq:c^hShPstar}, \eqref{eq:b^Sstar}, \eqref{eq:c^Sstar}, and \eqref{eq:201213.1} in Lemma~\ref{lem: bar to hat} that for all $\vt\hS_-=-1$ one has
\begin{align*}
b^{\Log(L)}+{}&(1+\Delta B^{\Log(L)})(X^{-2}_-\vt \tilde{c}^{S\star}\vt ^{\top}-2X^{-1}_-\vt b^{S\star})\\
&\qquad\qquad\qquad=b^{\Log(\hL),\hP}+(1+\Delta B^{\Log(\hL),\hP})(\vt \tilde{c}^{\hS,\hPstar}\vt ^{\top}-2\vt b^{\hS,\hPstar}).
\end{align*}
Since $\widehat{\Xi}_{\ha }$ in non-empty by Corollary~\ref{C:mainhP},  this shows $\Xi_a=\widehat{\Xi}_{\ha }/X_{-}$ is non-empty. 
Furthermore,
$$\frac{b^{\Log(L)}}{1+\Delta B^{\Log(L)}}=-\min_{\vt\in\R^{d}:\vt S_{-}=-1}\{\vt \tilde{c}^{S\star}\vt ^{\top}-2\vt b^{S\star}\}=2a b^{S\star}-a\tilde{c}^{S\star}a ^{\top},\qquad a\in\Xi_{a},$$
if and only if
$$\frac{b^{\Log(\hL),\hP}}{1+\Delta B^{\Log(\hL),\hP}}=-\min_{\vt\in\R^{d}:\vt \hS_{-}=-1}\{\vt \tilde{c}^{\hS,\hPstar}\vt ^{\top}-2\vt b^{\hS,\hPstar}\}=2\ha  b^{\hS,\hPstar}-\ha \tilde{c}^{\hS,\hPstar}\ha  ^{\top},\qquad 
\ha \in\hXi_{\ha }.$$

Next, setting $a=\ha /X_-$ for any $\ha \in\hXi_{\ha }$, we have, by Corollary \ref{cor: Exp equality}, that 
\begin{equation}
\Exp((-\ha \one_{\zu\tau ,T\zu})\sint \hS)X/X^{\tau}=\Exp((-\ha \one_{\zu\tau ,T\zu}/X_-)\sint S)=\Exp((-a\one_{\zu\tau ,T\zu})\sint S)\label{eq:expnc}
\end{equation}
for any $[0,T]$-valued stopping time $\tau$. Thanks to \eqref{eq:210421}, 
$-a\one_{\zu\tau ,T\zu}\Exp((-a\one_{\zu\tau ,T\zu})\sint S)_-\in \adm^{\tau}_{1}(S,\P)$ if and only if 
$-\ha \one_{\zu t ,T\zu}\Exp\big((-\ha \one_{\zu t ,T\zu})\sint \hS\big)_{-}\in\adm^{\tau}_{1}(\hS,\hP)$.

It follows that $L$ satisfies \ref{T:main.1} for the price process $S$ under $\P$ if and only if $\hL$ satisfies item~\ref{C:mainhP.1} of Corollary~\ref{C:mainhP}. Hence, $L$ satisfies \ref{T:main.1} if and only if it is the opportunity process for the price process $S$ under $\P$.

\ref{T:main.2}--\ref{T:main.3}: Combining $\hH=\frac{H}{X_T}$ and $L=\frac{\hZ}{X^2}\hL$ with \eqref{eq:expnc} and Bayes' formula, \eqref{eq:hV} and \eqref{eq:V} yield that
\begin{equation}\label{eq:V to hV}
V_{t}=\frac{1}{L_t}\E\big[\Exp\big((-\one_{\zu t,T\zu}a)\sint S\big)_{T}H\big|\F_{t}\big]
=X_t\frac{1}{\hL_t}\E^{\hP}\big[\Exp\big((-\one_{\zu t,T\zu}\ha )\sint \hS\big)_{T}\hH\big|\F_{t}\big]=X_t\hV_t.
\end{equation}
From \eqref{eq:c^hShPstar}, \eqref{eq:201212.1}, \eqref{eq:tchVhPstar}, \eqref{eq:c^Sstar}, \eqref{eq:c^SVstar}, \eqref{eq:c^Vstar}, and \eqref{eq:hatcvtobar} in Lemma~\ref{lem: bar to hat} one obtains 
\begin{equation}\label{eq:error to herror}
\begin{split} 
(1+\Delta B^{\Log(L)})X^{-2}_-(\vt \tilde{c}^{S\star}\vt ^{\top}-&2_-\vt \tilde{c}^{SV\star}+\tilde{c}^{V\star})\\
&{}=(1+\Delta B^{\Log(\hL),\hP})(\vt \tilde{c}^{\hS,\hP^\star}\vt ^{\top}-\vt c^{\hS\hV,\hP^\star}+\tilde{c}^{\hV,\hP^\star}),
\end{split}
\end{equation}
which shows $\Xi_{\xi}=\hXi_{\smash{\hxi}}$. From \eqref{eq:3.5} one has $\hvp(\hv ,\hH)=\vp(v,H)$. Earlier, we have established $\Xi_a=\hXi_{\ha}/X_-$ and in \eqref{eq:V to hV} also $\hV=V/X$. These observations and \eqref{eq:hphi} now yield
\begin{align*}
\vp(v,H)=\hvp(\hv ,\hH)&{}=\hxi +(\hV_{-}-\hvp (\hv ,\hH)\hS_{-})\ha \\
&{}=\hxi +(V_{-}-\hvp (\hv ,\hH)S_{-})\ha /X_-
{}=\xi +(V_{-}-\vp (v,H)S_{-})a,
\end{align*}
where $\xi$ and $a$ are arbitrary elements of $\Xi_a$ and $\Xi_\xi$, respectively. This proves \eqref{eq:mvhst}.

For the unconditional hedging error, we have from \eqref{eq:3}
\begin{align*}
\E[(\vp (v,H)S_{T}-H)^{2}]&=\E[X_T^2]\E^{\hP}[(\vp (v,H)\hS_{T}-\hH)^{2}]=\hZ_0\E^{\hP}[(\hvp(\hv ,\hH)\hS_{T}-\hH)^{2}]\\
&=\hZ_0\E^{\hP}\left[\hL_0(\hv -\hV_0)^2+\hL\big(\hxi\tilde{c}^{\hS,\hPstar}\hxi ^{\top}-2\hxi\tilde{c}^{\hS\hV,\hPstar}+\tilde{c}^{\hV,\hPstar}\big)\sint A_T\right],
\end{align*}
where the last equality follows from \eqref{eq:mvherrorhP}. Since $\big(\hxi\tilde{c}^{\hS,\hPstar}\hxi ^{\top}-2\hxi\tilde{c}^{\hS\hV,\hPstar}+\tilde{c}^{\hV,\hPstar}\big)\sint A$ is predictable, non-decreasing and starting at $0$ and hence locally of integrable variation, combining Rogers and Williams~\cite[Theorem~VI.21.1]{rogers.williams.94.vol2} with an application of the monotone convergence theorem and a localisation argument yields
\begin{align}
\E[(\vp (v,H)S_{T}-H)^{2}]&=\hZ_0\E^{\hP}\left[\hL_0(\hv -\hV_0)^2+{}^{\mathrm{p},\hP}\hL\big(\hxi\tilde{c}^{\hS,\hPstar}\hxi ^{\top}-2\hxi\tilde{c}^{\hS\hV,\hPstar}+\tilde{c}^{\hV,\hPstar}\big)\sint A_T\right]\label{eq:mvh:hat2}
\end{align}
for the predictable projection ${}^{\mathrm{p},\hP}(\hL)$ of $\hL$ under $\hP$. 

Because ${}^{\mathrm{p},\hP}(\hL)=\hL_{-}(1+\Delta B^{\Log(\hL),\hP})=L_-X_-^2\hZ_-^{-1}(1+\Delta B^{\Log(\hL),\hP})$, $\frac{\d\hP}{\d\P}=\frac{\hZ_T}{\hZ_0}$, and $\Xi_{\xi}=\hXi_{\smash{\hxi}}$, we obtain from \eqref{eq:error to herror} that \eqref{eq:mvh:hat2} equals
$$
\E\left[\hZ_T\left(L_0X_0^2\hZ_0^{-1}(\hv -\hV_0)^2
+L_-\hZ_-^{-1}(1+\Delta B^{\Log(L)})(\xi\tilde{c}^{S\star}\xi ^{\top}-2\xi\tilde{c}^{SV\star}+\tilde{c}^{V\star})\sint A_T\right)\right]
$$
and hence also
$\E\left[L_0(c-V_0)^2+L(\xi\tilde{c}^{S\star}\xi ^{\top}-2\xi\tilde{c}^{SV\star}+\tilde{c}^{V\star})\sint A_T\right]$
for arbitrary $\xi\in\Xi_\xi$. The latter follows by applying \cite[Theorem~VI.21.1]{rogers.williams.94.vol2} twice together with the monotone convergence theorem and a localisation argument, since 
$
(\xi\tilde{c}^{S\star}\xi ^{\top}-2\xi\tilde{c}^{SV\star}+\tilde{c}^{V\star})\sint A
$ 
is predictable and non-decreasing and hence locally of integrable variation. In the first application of \cite[Theorem~VI.21.1]{rogers.williams.94.vol2}, we use that $\hZ$ is a martingale under $\P$ so that ${}^{\mathrm{p}}(\hZ_T)=\hZ_-$ holds for the predictable projection and in the second that $^{\mathrm{p}}L=L_-(1+\Delta B^{\Log(L)})$. This shows \eqref{eq:mvherror}, which completes the proof.
\newpage
\addtocounter{appx}{1}
\section{Auxiliary statements for Theorem~\ref{T:main}.}\label{App:D}
In this appendix, $S=(X,Y)$ is locally square-integrable under $\P$, $X>0$, $X_->0$, $X^2_T\in L^1(\P)$, $\hP$ is defined by $\d\hP/\d\P=X_T^2/\E[X_T^2]$, $\hZ_t=\E[X_T^2|\F_t]$ for $0\leq t\leq T$, and $\hS = S/X$ is locally square-integrable under $\hP$.
\begin{lemma}\label{lem: bar to hat}
Let $(\hL ,\hV )$ and $(L,V)$ be semimartingales under $\hP$ and $\P$, respectively, such that $\hL _->0$, $L_->0$, 
\begin{align*}
\hb ^{\hS}   &{}=b^{\hS,\hP}+c^{\hS\Log(\hL )}+\int ylF^{\hS,\Log(\hL ),\hP}\big(\d(y,l)\big), \\
\hc ^{\hS}   &{}=c^{\hS}+\int yy^{\top }(1+l)F^{\hS,\Log(\hL ),\hP}\big(\d(y,l)\big), \\
\hc ^{\hV\hS}&{}=c^{\hV\hS}+\int x^{\top }z(1+y)F^{\hS,\Log(\hL),\hV,\hP}\big(\d(x,y,z)\big),\\
\hc ^{\hV}   &{}=c^{\hV}+\int z^2(1+y)F^{\Log(\hL),\hV,\hP}\big(\d(y,z)\big),
\end{align*}
and
\begin{align*}
\bar{b}^S   &{}= b^{S}+c^{S\Log(L)}+\int xyF^{S,\Log(L)}\big(\d(x,y)\big) \\ 
\bar{c}^S   &{}= c^{S}+\int xx^{\top }(1+y)F^{S,\Log(L)}\big(\d(x,y)\big) \\ 
\bar{c}^{VS}&{}= c^{VS}+\int x^{\top }z(1+y)F^{S,\Log(L),V}(\d(x,y,z)),\\
\bar{c}^{V} &{}= c^{V}+\int z^2(1+y)F^{\Log(L),V}\big(\d(y,z)\big)
\end{align*}
are well-defined. Then, if $L=\frac{\hL \hZ }{X^2}$ and $V=X\hV $, we have that
\begin{align}
b^{\Log(\hL),\hP}&{}= b^{\Log(L)} + 2 X_-^{-1}\bar{b}^{X} +X_-^{-2}\bar{c}^{X}; \label{eq: b Log(hL) hP} \\
\hb ^{\hS}    &{}= X_{-}^{-1}\left(\bar{b}^S+X_{-}^{-1}\bar{c}^{SX}-X_{-}^{-1}S_-(\bar{b}^X +X_{-}^{-1}\bar{c}^X )\right);
\label{eq: bhat bbar}\\
\hc ^{\hS}    &{}=X_{-}^{-2}\left(\bar{c}^S-X_{-}^{-1}(S_{-}\bar{c}^{XS}+\bar{c}^{SX}S_{-}^\top)+X_{-}^{-2}\bar{c}^XS_{-}S_{-}^\top  \right); \label{eq: chat cbar}\\
\hc ^{\hS\hV }&=X_{-}^{-2}\left(\bar{c}^{SV}-X_{-}^{-1}(S_{-}\bar{c}^{VX}+\bar{c}^{SX}V_{-})+X_{-}^{-2}\bar{c}^XV_{-}S_{-} \right); \label{eq:190326.2}\\
\hc ^{\hV}&=X_{-}^{-2}\left(\bar{c}^V-2\hV_{-}\bar{c}^{VX}+\hV_{-}^{2}\bar{c}^X  \right). \label{eq: cvhat cvbar}
\end{align}
Moreover, we have
\begin{equation}\label{eq:201213.1}
 \vt\hc ^{\hS}\vt^\top-2\vt\hb ^{\hS}+b^{\Log(\hL),\hP} = X_{-}^{-2}\vt\bar{c}^S\vt^\top -2X_{-}^{-1}\vt\bar{b}^S+b^{\Log(L)}
,\qquad \text{for all $\vt \hS_-=-1$},
\end{equation}
and 
\begin{align}\label{eq:hatcvtobar}
\vt\hc^{\hS}\vt^\top-2\vt\hc^{\hS\hV}+\hc^{\hV}=X^{-2}_-\left(\vt\bar{c}^S\vt^\top-2\vt\bar{c}^{SV}+\bar{c}^V \right),\qquad
\text{for all $\vt \hS_-=\hV$}.
\end{align}
\end{lemma}
\bpf
Observe that 
\begin{alignat}{2}
 \bar{b}^S    \sint A&{}= B^{S+[S,\Log(L)]}; \qquad &  \bar{c}^S \sint A&{}= B^{[S,S]+\left[ [S,S],\Log(L)\right]}; \label{eq: barb1barc1}\\
 \bar{c}^{VS} \sint A&{}= B^{[V,S]+\left[ [V,S],\Log(L)\right]};\qquad& \bar{c}^{V}\sint A&{}= B^{[V,V]+\left[ [V,V],\Log(L)\right]}.
\label{eq: barcv1barv1}
\end{alignat}
By Lemma~\ref{L:aux1} one has
$$ b^{\Log(\hL),\hP} = b^{\Log(\hL)+[\Log(\hL),\Log(\hZ)]} = b^{\Log(\hZ)+\Log(\hL)+[\Log(\hL),\Log(\hZ)]} 
= b^{\Log(\hL\hZ)}=b^{\Log(LX^2)},$$
where the second equality holds due to $\hZ$ being a local $\P$-martingale, the third follows from the Yor formula, and the last from the identity $\hL\hZ=LX^{2}$. A double application of the Yor formula yields
$$ \Log(LX^2) = \Log(L) +2\Log(X)+[\Log(X),\Log(X)]+[2\Log(X)+[\Log(X),\Log(X)],\Log(L)].$$
On matching individual terms to expressions in \eqref{eq: barb1barc1} one finally obtains \eqref{eq: b Log(hL) hP}.

Analogously,
\begin{alignat*}{2}
\hb ^{\hS}   \sint A&{}= B^{\hS+[\hS,\Log(\hL)],\hP};\qquad& \hc ^{\hS}\sint A&{}=B^{[\hS,\hS]+[ [\hS,\hS],\Log(\hL)],\hP};\\
\hc ^{\hV\hS}\sint A&{}= B^{[\hV,\hS]+[ [\hV,\hS],\Log(\hL)],\hP};\qquad&  \hc ^{\hV}\sint A&{}=B^{[\hV,\hV]+[ [\hV,\hV],\Log(\hL)],\hP},
\end{alignat*}
whereby, in view of $\hL\hZ=LX^{2}$, Lemma \ref{L:aux1} yields
\begin{alignat}{2}
\hb ^{\hS}   \sint A&{}= B^{\hS+[\hS,\Log(LX^{2})]};\qquad&\hc ^{\hS}\sint A&{}= B^{[\hS,\hS]+[[\hS,\hS],\Log(LX^{2})]};   
\label{eq: hatb2hatc2}\\
\hc ^{\hV\hS}\sint A&{}= B^{[\hV,\hS]+[ [\hV,\hS],\Log(LX^2)\big]};\qquad&\hc ^{\hV}\sint A&{}= B^{[\hV,\hV]+[ [\hV,\hV],\Log(LX^2)]}.    \label{eq: hatcvs2hatcvv2}
\end{alignat}
From  Yor's formula, 
\begin{equation} \label{eq: Yor}
\Log(X)+\Log(X^{-1})+[\Log(X),\Log(X^{-1})]=0,
\end{equation}
and integration by parts one obtains
\begin{align}
\hS -\hS_0    &{}=X_{-}^{-1}\sint S +S_-\sint X^{-1} +[S,X^{-1}]\nonumber\\
              &{}=X_{-}^{-1}\sint (S +S_-\sint \Log(X^{-1}) +[S,\Log(X^{-1})])
              {}=X_{-}^{-1}\sint (\tilde{S}+[\tilde{S},\Log(X^{-1})]),\label{eq: hatS tildeS}
\end{align}
with $\tilde{S}=S-S_{-}\sint \Log(X)$. 
Likewise, one obtains
\begin{equation}\label{eq: hatV tildeV}
\hV=\hV_{0}+X_{-}^{-1}\sint \big(\tilde{V}+[\tilde{V},\Log(X^{-1})]\big)
\end{equation}
with $\tilde{V} =V-V_{-}\sint \Log(X)$.
From \eqref{eq: hatb2hatc2}, \eqref{eq: hatS tildeS}, and Lemma \ref{lem: aux2} with $\tilde{W}=X_{-}^{-1}\sint \tilde{S}$ one 
obtains
\begin{align}
\hb ^{\hS}\sint A
&{}=B^{\tilde{W}+[\tilde{W},\Log(X)]+[\tilde{W}+[\tilde{W},\Log(X)],\Log(L)]}.\label{eq: hatb3}
\end{align}
In view of 
$\tilde{W}+[\tilde{W},\Log(X)]=X_{-}^{-1}\sint\left(S-S_-\sint \Log(X)+[S,\Log(X)]-S_-\sint [\Log(X),\Log(X)] \right)$,
a combination of \eqref{eq: barb1barc1} and \eqref{eq: hatb3} yields \eqref{eq: bhat bbar}.

Similarly, from \eqref{eq: hatS tildeS} and \eqref{eq: hatV tildeV}, one has
$$
[ \hV ,\hS] =X_{-}^{-2}\sint \left([\tilde{V},\tilde{S}]+2[ \tilde{V},[\tilde{S},\Log(X^{-1})] ]+[ [\tilde{V},\Log(X^{-1})],[\tilde{S},\Log(X^{-1})] ]\right),$$
which one can rephrase as 
\begin{equation}\label{eq: [hatS,hatS]}
[ \hV ,\hS] =\tilde{W}+[\tilde{W},\Log(X^{-1})],
\end{equation}
with $\tilde{W}=X_{-}^{-2}\sint([\tilde{V},\tilde{S}]+[ [\tilde{V},\tilde{S}],\Log(X^{-1})])$. 

To obtain \eqref{eq:190326.2}, apply Lemma~\ref{lem: aux2} to \eqref{eq: hatcvs2hatcvv2}, making use of \eqref{eq: [hatS,hatS]}. This then yields
\begin{align*}
\hc ^{\hV \hS}\sint A &{}=B^{\tilde{W}+[\tilde{W},\Log(X)]+[\tilde{W}+[\tilde{W},\Log(X)],\Log(L)]} 
=X_{-}^{-2}\sint B^{[ \tilde{V},\tilde{S}]+[ [\tilde{V},\tilde{S}],\Log(L)]}
\end{align*}
on observing that thanks to \eqref{eq: Yor} one has $\tilde{W}+[\tilde{W},\Log(X)] =X_{-}^{-2}\sint [ \tilde{V},\tilde{S}]$. 
Formula~\eqref{eq:190326.2} now follows by substituting for $\tilde{S},\tilde{V}$ and matching terms in \eqref{eq: barb1barc1}--\eqref{eq: barcv1barv1}.
Formulae \eqref{eq: chat cbar} and \eqref{eq: cvhat cvbar} are obtained analogously.
 
For $\vt\hS_{-}=-1$, one has $\vt S_{-}=-X_-$ and therefore by \eqref{eq: bhat bbar} and \eqref{eq: chat cbar}
\begin{align*}
2\vt\hb ^{\hS}-\vt\hc ^{\hS}\vt^\top &{}= 2X_{-}^{-1}\vt\bar{b}^S-X_{-}^{-2}\vt\bar{c}^S\vt^\top
+2\bar{b}^X X_{-}^{-1} +\bar{c}^X X_{-}^{-2}.
\end{align*}
Together with \eqref{eq: b Log(hL) hP} this yields \eqref{eq:201213.1}.

Since $V=X\hV$, we have that $\vt S_-=V_-$ if and only if $\vt \hS_-=\hV_-$. Therefore, for $\R^d$-valued $\vt$ with $\vt S_-=V_-$ and hence $\vt \hS_-=\hV_-$, one obtains  from \eqref{eq: chat cbar} and \eqref{eq:190326.2} after simplifications
$$\vt\hc^{\hS}\vt^\top-2\vt\hc^{\hS\hV}=X^{-2}_-\left(\vt\bar{c}^S\vt^\top-2\vt\bar{c}^{SV}+2\hV_-\bar{c}^{VX}-\hV_-^2\bar{c}^X\right).$$
Together with \eqref{eq: cvhat cvbar} this yields \eqref{eq:hatcvtobar}.
\ep

\begin{lemma}\label{L:aux1} 
For any $\hL$, $\hZ$ as above and any semimartingale $\hW$, the following are equivalent.
\begin{enumerate}[(1)]
\item $\hW +[\hW ,\Log(\hL )]$ is ${\hP}$-special.
\item $\hW +[\hW ,\Log(\hL \hZ )]$ is ${\P}$-special.
\end{enumerate}
If one (hence both) conditions hold, then
$ B^{\hW +[\hW ,\Log(\hL )],\hP}=B^{\hW +[\hW ,\Log(\hL \hZ )]}$. 
\end{lemma}

\bpf
By Girsanov's theorem for equivalent measures, $\hW +[\hW ,\Log(\hL )]$ is ${\hP}$-special if and only if
$ \hW +[\hW ,\Log(\hL )]+[\hW +[\hW ,\Log(\hL )],\Log(\hZ )] \text{ is $\P$-special}$, 
in which case the $\hP$ compensator of the former equals the $\P$-compensator of the latter. A simple manipulation together with Yor's formula now yield
\begin{align*}
\hW +[\hW ,\Log(\hL )]+[\hW +[\hW ,\Log(\hL )],\Log(\hZ )] 
&{}=\hW +[\hW ,\Log(\hL )+\Log(\hZ )+[\Log(\hL ),\Log(\hZ )] ] \\
&{}=\hW +[\hW ,\Log(\hL \hZ )],
\end{align*}
which completes the proof.
\ep

\begin{lemma}
\label{lem: aux2}For $X, X_- >0$, $L, L_- >0$ and any semimartingale $\tilde{W}$ we have
\begin{align*}
\tilde{W}+{}[\tilde{W},\Log(X^{-1})]+{}&[\tilde{W}+[\tilde{W},\Log(X^{-1})],\Log(LX^{2})] =\tilde{W}+[\tilde{W},\Log(LX)]\\
={}&\tilde{W}+[\tilde{W},\Log(X)]+[\tilde{W}+[\tilde{W},\Log(X)],\Log(L)].
\end{align*}
\end{lemma}
\bpf By direct calculation one obtains
\begin{align*}
\tilde{W}+{}&[\tilde{W},\Log(X^{-1})]+[\tilde{W}+[\tilde{W},\Log(X^{-1})],\Log(LX^{2})] \\
={}&\tilde{W}+[\tilde{W},\Log(X^{-1})+\Log(LX^{2})+[\Log(X^{-1}),\Log(LX^{2})] ] 
={}\tilde{W}+[\tilde{W},\Log(LX)] \\
={}&\tilde{W}+[ \tilde{W},\Log(X)+\Log(L)+[\Log(X),\Log(L)] ] 
={}\tilde{W}+[\tilde{W},\Log(X)]+[\tilde{W}+[\tilde{W},\Log(X)],\Log(L)].
\end{align*}
Here the second and third equality follow from the Yor formula.
\ep

\begin{lemma} \label{lem: sf Exp}
Suppose $\alpha\in L(S)$ and $\alpha S_{-}=1$. Then $\vt = \alpha\Exp(\alpha \sint S)_{-}$ is a self-financing strategy with $\vt S=\Exp(\alpha \sint S)$. Furthermore, for $\hat{\alpha}=\alpha X_-$ one has $\hat\alpha\in L(\hS)$ and
\begin{equation}\label{eq: Exp equality} 
\Exp(\alpha \sint S)=\Exp(\hat{\alpha} \sint \hS)X/X_0.
\end{equation}
\end{lemma}
\bpf
By assumption, we have $\vt_0 S_0=1$ and $\vt S_{-}=\Exp(\alpha \sint S)_{-}$. From the properties of stochastic exponential one obtains $\vt S_{-} = 1 + \vt\sint S_-$, hence $\vt$ is self-financing by Proposition~\ref{P:sf}. This also yields $\vt S = 1 + \vt\sint S=\Exp(\alpha \sint S)$, which completes the proof of the first assertion. 

Since $\alpha\in L(S)$ and $\alpha S_-=1$, we have $\hat\alpha\in L((1/X_-)\sint S)\cap L([S,1/X])\cap L(S_-\sint (1/X))$. This yields $\hat\alpha\in L(S/X)$ through integration by parts. By Proposition~\ref{P:sf}\ref{P:sf.iii} we have $\vt\in L(\hS)$ and $\vt S/X=\vt_0 S_0/X_0 + \vt\sint \hS$. 
Multiplying by $X_0$ and substituting $\vt S=\Exp(\alpha\sint S)$ gives 
$$
X_0 \Exp(\alpha\sint S)/X = 1 +\alpha X_0 \Exp(\alpha\sint S)_{-}\sint\hS
=1 + (X_0\Exp(\alpha\sint S)_{-}/X_{-})\sint(\hat\alpha \sint\hS),
$$ 
yielding $X_0\Exp(\alpha\sint S)/X = \Exp(\hat\alpha \sint\hS)$. 
\ep
\begin{corollary}\label{cor: Exp equality}
Let $\tau$ be a stopping time such that $\one_{\zu \tau,T\zu}\hat{\alpha}\in L(\hS)$ and $\hat{\alpha} \hS_{-}=1$ on $\zu \tau,T\zu$. Then 
$$\Exp(\one_{\zu \tau,T\zu}\hat{\alpha}/X_{-}\sint S)=\Exp(\one_{\zu \tau,T\zu}\hat{\alpha} \sint \hS)X/X^\tau.$$
\end{corollary}

\bpf
On $\auf 0,\tau\zu$ let $\hat{\alpha} = [1,0,\ldots,0]$ and denote $\alpha = \hat{\alpha}/X_{-}$. Then $\alpha S_{-}=1$ on $[0,T]$ and \eqref{eq: Exp equality} applies. By the Yor formula $\Exp(Z)=\Exp(\one_{\auf 0,\tau\zu}\sint Z)\Exp(\one_{\zu \tau,T\zu}\sint Z)$ for any semimartingale $Z$, while $\Exp(\one_{\auf 0,\tau\zu}\alpha\sint S)=X^\tau/X_0$. Since $X>0$, the statement follows.
\ep
\newpage
\addtocounter{appx}{1}
\section{Glossary of mathematical notation and terminology.}\label{App:E}\leavevmode
{} \\[1ex]
{\small
\begin{glossb}
  \item[$A$]    activity process (predictable, increasing, and integrable)
  \item[$B$, $b$]    cumulative drift and its rate relative to $A$
  \item[$C$, $c$]    continuous quadratic variation and its rate relative to $A$
	\item[$\Exp$, $\Log$] stochastic exponential, stochastic logarithm
	\item[$L$] the opportunity process
	\item[$L(X)$] the set of predictable $X$-integrable processes
	\item[$L^\infty(\P)$, $L^2(\P)$] the set of bounded (resp., square-integrable) random variables
	\item[$\cM_2(S,\P)$] the set of deflators for $S$ under $\P$
	\item[$\nu$, $F$]  predictable compensator of the jump measure and its rate relative to $A$
	\item[$\Null$] null space of a matrix
	\item[$\NullStrategy$] null strategies (trading strategies whose wealth is identically 0)
	\item[$\P,\P^\star\!,\hP,\hP^\star\!,\tilde \P$] probability measures
	\item[${}^\mathrm{p}L$, ${}^{\mathrm{p},\hP}\hL$] predictable projection of $L$ under $\P$, resp., that of  $\hL$ under $\hP$
	\item[$\Ran$] range of a matrix
	\item[$\rho$, $\sigma$, $\tau$] stopping times
	\item[$T$] terminal date
	\item[$S$, $\hS$] price processes
	\item[$\bot$] orthogonal complement
	\item[$\oplus,\oplus^\bot$] direct sum (resp., direct orthogonal sum) of two subspaces
	\item[$\top$] matrix transpose
	\item[$\Theta$, $\adm$] tame and admissible trading strategies
	\item[$\vt \sint S_t$] stochastic integral $\int_{(0,t]}\vt_u \d S_u$
	\item[$\zu\tau,T\zu$, $\auf\tau,T\zu$] stochastic intervals
	\item[$\ones$] column vector of ones
	\item[$\one$] indicator function
	\item[localization]    the procedure of stopping a process at a  sequence of stopping times increasing to $\infty$
  \item[locally square-integrable]    process that becomes square-integrable after localization
	\item[local martingale]    process that becomes a martingale after localization
	\item[truncation function]          function used to split jumps into  small vs. large
\end{glossb}
}
\bigskip
\end{APPENDICES}

\section*{Acknowledgments.}
We wish to thank two anonymous referees and an associate editor for their helpful comments. 

\newpage 

\def\shorturl#1{\href{http://#1}{\nolinkurl{#1}}}
\def\MR#1{\href{http://www.ams.org/mathscinet-getitem?mr=#1}{MR#1}}
\def\ARXIV#1{\href{https://arxiv.org/abs/#1}{arXiv:#1}}
\def\DOI#1{\href{https://doi.org/#1}{doi:#1}}

\end{document}